\documentclass[11pt,a4paper]{amsart}

\usepackage{amsmath,amssymb,amsthm, amsfonts,enumerate,mathtools,pinlabel,float,stmaryrd}
\usepackage[mathscr]{eucal}
\usepackage{amsbsy}
\usepackage{graphicx}
\usepackage{caption}
\usepackage{subcaption}
\usepackage{cite}
\usepackage{array}
\usepackage{xcolor}
\usepackage{calc}
\usepackage{comment}
\usepackage{soul}
\usepackage[myheadings]{fullpage}
\usepackage{hyperref}
\hypersetup{
    colorlinks=true,   
    linkcolor=blue,
    citecolor=blue,
}

\newtheorem{theorem*}{Theorem}
\newtheorem{cor*}{Corollary}
\newtheorem{prop*}{Proposition}

\newtheorem{theorem}{Theorem}[section]
\newtheorem{cor}{Corollary}[theorem]
\newtheorem{prop}[theorem]{Proposition}
\newtheorem{lem}[theorem]{Lemma}
\newtheorem{algo}[theorem]{Algorithm}

\theoremstyle{definition}
\newtheorem{defn}[theorem]{Definition}
\newtheorem{exmp}[theorem]{Example}
\newtheorem{rem}[theorem]{Remark}

\DeclareMathOperator{\lcm}{lcm}
\DeclareMathOperator{\F}{\mathcal{F}}

\DeclareMathOperator{\map}{Mod}
\DeclareMathOperator{\homeo}{Homeo^{+}}
\DeclareMathOperator{\sgn}{sgn}
\everymath{\displaystyle}

\makeatletter
\@namedef{subjclassname@2020}{\textup{2020} Mathematics Subject Classification}
\makeatother

\begin{document}

\title{General primitivity in the mapping class group}

\author{Pankaj Kapari}
\address{Department of Mathematics\\
Indian Institute of Science Education and Research Bhopal\\
Bhopal Bypass Road, Bhauri \\
Bhopal 462 066, Madhya Pradesh\\
India}
\email{pankajkapri02@gmail.com}
\urladdr{https://sites.google.com/view/pankajkapdi/home}

\author{Kashyap Rajeevsarathy}
\address{Department of Mathematics\\
Indian Institute of Science Education and Research Bhopal\\
Bhopal Bypass Road, Bhauri \\
Bhopal 462 066, Madhya Pradesh\\
India}
\email{kashyap@iiserb.ac.in}
\urladdr{https://home.iiserb.ac.in/$_{\widetilde{\phantom{n}}}$kashyap/}



\begin{abstract}
For $g\geq 2$, let $\mathrm{Mod}(S_g)$ be the mapping class group of the closed orientable surface $S_g$ of genus $g$. In this paper, we obtain necessary and sufficient conditions under which a given pseudo-periodic mapping class can be a root of another up to conjugacy. Using this characterization, the canonical decomposition of (non-periodic) mapping classes, and some known algorithms, we give an algorithm for determining the conjugacy classes of roots of arbitrary mapping classes. Furthermore, we derive realizable bounds on the degrees of roots of pseudo-periodic mapping classes in $\mathrm{Mod}(S_g)$, the Torelli group, the level-$m$ subgroup of $\mathrm{Mod}(S_g)$, and the commutator subgroup of $\mathrm{Mod}(S_2)$. In particular, we show that the highest possible (realizable) degree of a root of a pseudo-periodic mapping class $F$ is $3q(F)(g+1)(g+2)$, where $q(F)$ is a unique positive integer associated with the conjugacy class of $F$. Moreover, this bound is realized by a root of a power of a Dehn twist about a separating curve of genus $[g/2]$ in $S_g$, where $g\equiv 0,9 \pmod{12}$. Finally, for $g\geq 3$, we show that any pseudo-periodic mapping class having a nontrivial periodic component that is not the hyperelliptic involution, normally generates $\mathrm{Mod}(S_g)$. Consequently, we establish that $\mathrm{Mod}(S_g)$ is normally generated by a root of bounding pair map or a root of a nontrivial power of a Dehn twist. 
\end{abstract}

\maketitle

\section{Introduction}
\label{sec:intro}
For $g\geq 2$, let $\map(S_g)$ be the \textit{mapping class group} of the closed orientable surface $S_g$ of genus $g$. A nontrivial $G\in \map(S_g)$ is said to be a \textit{root of an $F\in \map(S_g)$ of degree $n$} if there exists a least integer $n>1$ such that $G^n=F$ and $|G|=n|F|$. If $F$ does not have roots of degree $n$ for any $n$, then it is said to be \textit{primitive}. A natural question in this context is whether one can determine if an arbitrary $F \in \map(S_g)$ is primitive, and compute the conjugacy classes of the roots of $F$ when it is non-primitive. We call this the \textit{general primitivity problem} in $\map(S_g)$.

Fehrenbach and Los \cite{los,los1} have developed an algorithm that computes the roots of pseudo-Anosov mapping classes on once-punctured surfaces. Furthermore, there is an algorithm known as the \textit{flipper algorithm} due to Bell \cite{flipper,flipper1} that can determine the primitivity of pseudo-Anosov mapping classes. This algorithm works by finding the Agol veering triangulation \cite{agol} of the mapping torus. For non-primitive mapping classes a root can be extracted from a subset of the triangulation of this bundle. For the case of a single Dehn twist (or its power) and a product of commuting Dehn twists, the answer to this question is well known\cite{margalit, rajeevsarathy3, monden, rajeevsarathy4, rajeevsarathy5, rajeevsarathy6}. More recently, Dhanwani and Rajeevsarathy \cite{dhanwani1} have given equivalent conditions for the primitivity of periodic mapping classes. In view of the canonical decomposition of non-periodic mapping classes, a key step in addressing the general primitivity problem in $\map(S_g)$ is to determine the primitivity of an arbitrary \textit{pseudo-periodic} mapping class (i.e. an infinite order reducible mapping class that has only periodic components in its canonical decomposition).

In this paper, to begin with, we classify the conjugacy classes of roots of \textit{multitwists} (i.e. a product of powers of commuting Dehn twists) in Section~\ref{sec:gen_prim}. This is a complete generalization of the theory developed in \cite{rajeevsarathy6} (see Proposition~\ref{thm:main_theorem}). Furthermore, we derive equivalent conditions under which a given pseudo-periodic mapping class can be a root of another (see Theorem \ref{thm:pp_primitive}) up to conjugacy. For obtaining these results, we apply Thurston's orbifold theory \cite[Chapter 13]{thurston1} and theory of pseudo-periodic mapping classes developed in \cite{rajeevsarathy6, matsumoto}.

There are well known polynomial-time algorithms \cite{curver,margalit1} for determining the Nielsen-Thurston type of a mapping class given as a product of Lickorish Dehn twists. However, for this paper we will choose \textit{the curver algorithm} from~\cite{curver1} as this has already been implemented~\cite{curver} in the Python programming language. Further, this algorithm can compute the canonical reduction systems of arbitrary infinite order reducible mapping classes and the signature of the quotient orbifold of periodic mapping classes. This algorithm works by finding the shortest path between any two points of curve complex \cite{curver1}. By using our result on the primitivity of pseudo-periodic mapping classes, and the curver and flipper algorithms, we give an algorithm for solving the general primitivity problem in $\map(S_g)$.

\begin{theorem*}
For $g\geq 2$, there exists an algorithm that determines the primitivity of an arbitrary $F \in \map(S_g)$, and also computes the conjugacy classes of its roots, when $F$ is non-primitive.
\end{theorem*}

Let $T_c$ denote the left-handed Dehn twist about a simple closed curve $c$ in $S_g$. For an arbitrary pseudo-periodic $F$, there exists a unique positive integer $q(F)$ associated to its conjugacy class (see Proposition~\ref{thm:main_theorem}). For a pseudo-periodic mapping class $F$ that is a root of a multitwist $G = T_{c_1}^{q_1}T_{c_2}^{q_2}\cdots T_{c_m}^{q_m}$, this positive integer $q(F)$ will be defined as $\min\{|q_1|,|q_2|,\ldots,|q_m|\}$. In Section~\ref{sec:bounds}, as an application of these results, we obtain realizable bounds on the degrees of roots of pseudo-periodic mapping classes in $\map(S_g)$. A simple closed curve bounding a subsurface of least genus $g'$ in $S_g$ is called a \textit{separating curve of genus $g'$}.
 
\begin{prop*}
\label{prop:general_bound_i}
For $g\geq 2$, let $F\in \map(S_g)$ be a pseudo-periodic mapping class. Then the highest degree of a root of $F$ is $3q(F)(g+1)(g+2)$. Moreover, this bound is realized only when $F = T_c^{q}$, where $c$ is a separating curve of genus $[g/2]$ with $g\equiv 0,9 \pmod{12}$ and $q$ is a positive integer. In particular, $3(g+1)(g+2)$ is the highest degree of a root of a Dehn twist about a separating curve.
\end{prop*}

\noindent The following application characterizes the periodic components arising in the canonical decomposition of bound-realizing pseudo-periodic mapping classes.
 
\begin{cor*}
The upper bound-realizing roots of pseudo-periodic mapping classes in $\map(S_g)$ are primitive pseudo-periodic mapping classes that decompose canonically into irreducible periodic mapping classes whose Nielsen representatives have at least one fixed point.
\end{cor*}

\noindent We also obtain similar bounds for the pseudo-periodic mapping classes in the Torelli group $\mathcal{I}(S_g)$, the level-$m$ subgroup $\map(S_g)[m]$, and the commutator subgroup of $\map(S_2)$. A collection of nonseparating curves $C = \{c_1,c_2\}$ bounding a subsurface of least genus $g'$ in $S_g$ is called a \textit{bounding pair of genus $g'$} and the mapping class $T_{c_1}T_{c_2}^{-1}$ associated with $C$ is called a \textit{bounding pair map}. In Proposition \ref{prop:general_bound_i}, we have seen that a root of a Dehn twist about a separating curve realizes the highest degree in $\map(S_g)$. Since the powers of a Dehn twist about a separating curve lies in $\map(S_g)[m]$ for every $m$, in the following proposition we will assume that there are no separating curves in the canonical reduction system.

\begin{prop*}
\label{prop:torelli_levelm_i}
For $g \geq 3$ and $m \geq 2$, let $F\in \map(S_g)[m]$ be a pseudo-periodic mapping class. Suppose that the canonical reduction system for $F$ does not contain any separating curves and that $F$ has a root in $\map(S_g)$ of degree $n$.
\begin{enumerate}[(i)]
\item If $F\in \mathcal{I}(S_g)$, then $2 \leq n \leq q(F)g(g-2)$. The upper bound for $n$ is realized only when $F = (T_{c_1}T_{c_2}^{-1})^q$, where $\{c_1,c_2\}$ is a bounding pair of genus $(g/2)-1$ with $2\mid g$ and $q$ is a positive integer.
\item If $F\in \map(S_g)[m]\setminus \mathcal{I}(S_g)$, then $m \leq n \leq 3q(F)g(g-2)$. The upper bound for $n$ is realized when $F = (T_{c_1}T_{c_2}^k)^q$, where $\{c_1,c_2\}$ a bounding pair of genus $(g/2)-1$ with $g \equiv -4 \pmod{24}$, $q$ is a positive integer, and $k=(1/4)g(g-2) - 1$.
\end{enumerate}
In particular, $g(g-2)$ is the highest realizable degree of a root of a bounding pair map.
\end{prop*}

\noindent For $m\geq 3$, it is known \cite[Corollary 1.8]{ivanov} that any pseudo-periodic mapping class in $\map(S_g)[m]$ is a multitwist. By applying Proposition \ref{prop:torelli_levelm_i}, we obtain the following.

\begin{cor*}
For $g,m\geq 3$, let $F\in \map(S_g)[m]\setminus \mathcal{I}(S_g)$ be a pseudo-periodic mapping class. Suppose that the canonical reduction system for $F$ does not contain any separating curves and $F$ is not a power of a multitwist in $\map(S_g)$. If $m>3g(g-2)$, then $F$ is primitive.
\end{cor*}

\noindent As another application of our theory, we have the following. The commutator subgroup of $\map(S_2)$ will be denoted by $\map(S_2)^{(1)}$.

\begin{prop*}
Let $F\in \map(S_2)^{(1)}$ be a pseudo-periodic mapping class having a root in $\map(S_2)$ of degree $n$. Then we have
\begin{center}
$2q(F)\leq n\leq 12q(F)$.
\end{center}
Furthermore, the upper and lower bounds are realized when $F = T_c^q$, where $q$ is a positive integer such that $5$ divides $q$ and $c$ is a separating curve in $S_2$.
\end{prop*}

Recently, Margalit and Lanier \cite[Theorem 1.1]{normal} have proved that when $g\geq3$, any periodic mapping class, which is not a hyperelliptic involution, normally generates $\map(S_g)$. In Section~\ref{sec:normal_gen}, we generalize this result to pseudo-periodic mapping classes by applying the theory developed in \cite{rajeevsarathy1} and the \textit{well-suited curve criterion} from \cite{normal} (see Proposition \ref{prop:pp_normal}). Keeping in mind that the nontrivial powers of Dehn twists and bounding pair maps do not normally generate $\map(S_g)$ (as they act trivially on $H_1(S_g;\mathbb{Z}_k)$ for $k > 1$), we show the following as a final application.

\begin{prop*}
For $g \geq 3$, let $F \in \map(S_g)$ such that either $F$ is a nontrivial power of a Dehn twist or $F$ is a bounding pair map. Then there exists a root of $F$ that normally generates $\map(S_g)$.
\end{prop*}

\section{Preliminaries}
This section will introduce some basic notions pertaining to cyclic actions on surfaces, orbifold theory, and pseudo-periodic mapping class that are relevant to this paper. Throughout this paper we shall assume that $g\geq 2$.

\subsection{Periodic mapping classes} Let $F\in \map(S_g)$ be a periodic mapping class of order $n$. The Nielsen-Kerckhoff theorem \cite{kerckhoff} asserts that $F$ is represented by an orientation-preserving homeomorphism $\mathcal{F}$ of $S_g$ of the same order which induces a $C_n$-action on $S_g$, where $C_n=\langle \mathcal{F} \rangle$. The orbit space $\mathcal{O}_F:=S_g/C_n$ is the \textit{quotient orbifold} \cite[Chapter 13]{thurston1} associated to $F$ which is homeomorphic to $S_{g_0}$, where $g_0$ is the \textit{orbifold genus} of $\mathcal{O}_F$.

The $C_n$-action induces a branched covering $p:S_g \to \mathcal{O}_F$ with $k$ branch points (or \textit{cone points}) $x_1,x_2, \dots, x_k$ in the quotient orbifold $\mathcal{O}_F$ of orders $n_1,n_2,\dots, n_k$, respectively. The \textit{order} of a cone point $x_i$ is the order of the stabilizer subgroup of any point in the preimage of $x_i$. From orbifold covering space theory, the branch covering $p:S_g\to \mathcal{O}_F$ corresponds to an exact sequence 
\begin{equation*}
1 \longrightarrow \pi_1(S_g) \xrightarrow{p_*} \pi_1^{orb}(\mathcal{O}_F) \xrightarrow{\phi} C_n \longrightarrow 1.
\end{equation*}

\noindent Moreover, $\pi_1^{orb}(\mathcal{O}_F)$ is a Fuchsian group~\cite{katok,macbeath} that has following presentation:
\begin{equation*}
\label{eqn:fuchsian_present}
\langle \alpha_1,\beta_1,\alpha_2,\beta_2,\dots,\alpha_{g_0},\beta_{g_0},\gamma_1,\gamma_2\dots,\gamma_k\mid \gamma_1^{n_1}=\gamma_2^{n_2}=\dots=\gamma_k^{n_k}=\prod_{i=1}^k\gamma_i\prod_{i=1}^{g_0} [\alpha_i,\beta_i]=1 \rangle.
\end{equation*}
The homomorphism $\phi:\pi_1^{orb}\to C_n$ (classically known as \textit{surface kernel}) is an order-preserving epimorphism given by $\phi(\gamma_i)=F^{(n/n_i)d_i}$, where $d_i$ is a positive integer with $\gcd(d_i,n_i) =1$, for $1\leq i\leq k$.

The tuple $(g_0;n_1,n_2,\dots,n_k)$ is called the \textit{signature} of the quotient orbifold $\mathcal{O}_F$ and will be denoted by $\Gamma(\mathcal{O}_F)$. Each cone point $x_i$ of order $n_i$ in $\mathcal{O}_F$ lifts under $p$ to an orbit of size $n/n_i$ on $S_g$ and the \textit{local rotation} induced by $C_n$-action in this orbit is given by $2 \pi d_i^{-1}/n_i$, where $\gcd(d_i,n_i)=1$ and $d_id_i^{-1}\equiv 1 \pmod{n_i}$. Thus, a cyclic action along with the structure of its associated quotient orbifold can be compactly encoded as a tuple of integers.

\begin{defn}
\label{defn:cyclic_data}
A \textit{cyclic data set of degree $n$} is a tuple of the form 
\begin{center}
$D=(n,g_0;(d_1,n_1),(d_2,n_2),\dots,(d_k,n_k)),$
\end{center}
where $n\geq2$, $g_0\geq0$, $0\leq r \leq n-1$ with the following conditions.
\begin{enumerate}[(i)]
\item $n_i\geq2$, $d_i\geq 1$, $n_i\mid n$, $\gcd(d_i,n_i)=1$, $\text{ for all }i$.
\item $\lcm(n_1,n_2,\dots,\widehat{n_i},\dots,n_k)=\lcm(n_1,n_2,\dots,n_k)$, $\text{ for all }i$.
\item If $g_0=0$, then $\lcm(n_1,n_2,\dots,n_k)=n$.
\item $\sum_{i=1}^{k}\frac{n}{n_i}d_i\equiv 0 \pmod n $.
\item $\frac{2g-2}{n}=2g_0-2+\sum_{i=1}^{k}\left (1-\frac{1}{n_i}\right)$. \hfill \text{(Riemann-Hurwitz)}
\end{enumerate}
The number $g$ determined by the Riemann-Hurwitz equation is the \textit{genus} of the data set.
\end{defn}

The significance of the cyclic data set is given in the following proposition due to Nielsen \cite{nielsen} (see also~\cite[Theorem 3.9]{rajeevsarathy6}).
\begin{prop}
For $g\geq 2$, cyclic data sets of degree $n$ and genus $g$ are in one-to-one correspondence with the conjugacy classes of periodic mapping classes (or cyclic actions) of order $n$ in $\map(S_g)$.
\end{prop}

\begin{figure}[ht]
\tiny
\begin{subfigure}{\textwidth}
\labellist
\pinlabel $\pi$ at 165, 120
\pinlabel $\F$ at 153, 142
\endlabellist
\centering
\vspace{1mm}
\includegraphics[scale=0.6]{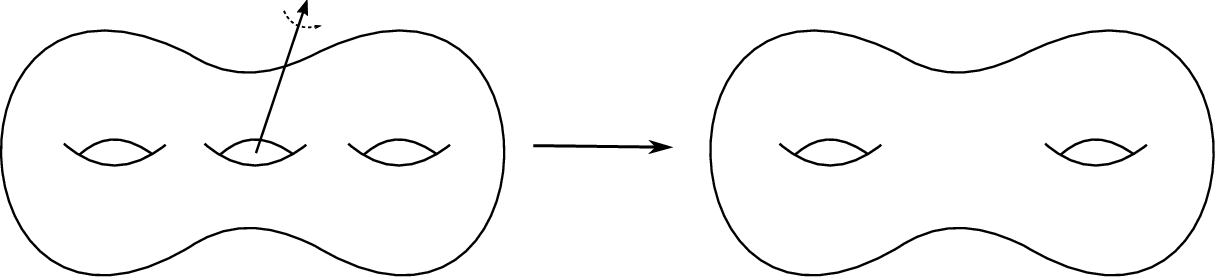}
\caption{A $\langle \F \rangle$-action on $S_3$ with data $D_{F}=(2,2;-)$.}
\end{subfigure}

\begin{subfigure}{\textwidth}
\labellist
\pinlabel $\F$ at 303, 60
\pinlabel $\pi$ at 290, 45
\endlabellist
\centering
\vspace{0.5cm}
\includegraphics[scale=0.6]{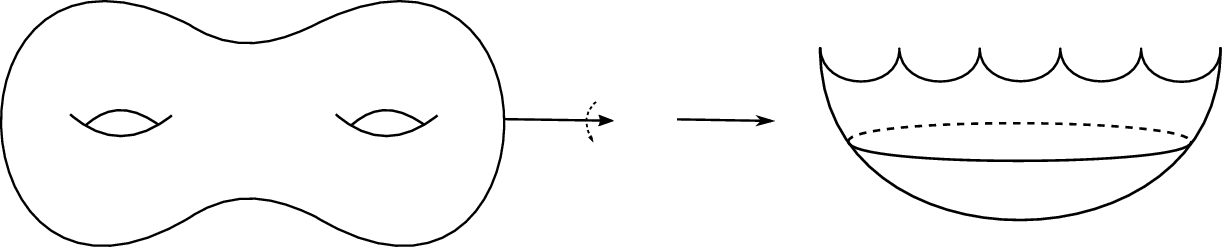}
\caption{A $\langle \F \rangle$-action on $S_2$ with data $D_{F}=(2,0;(1,2),(1,2),(1,2),(1,2),(1,2),(1,2))$.}
\end{subfigure}
\caption{Some illustration of cyclic data set.}
\label{fig1}
\end{figure}

\noindent The natural numbers $n$, $r$, $g$ and $g_0$ associated to a data set $D$ will be denoted by $n(D)$, $r(D)$, $g(D)$, and $g_0(D)$, respectively. From here on, a periodic mapping class $F$ (and its associated cyclic action) up to conjugacy will be represented by its corresponding data set $D_F$. Some cyclic actions on surfaces along with their corresponding data sets are shown in Figure \ref{fig1}.

\subsection{Pseudo-periodic mapping classes}
The Nielsen-Thurston classification \cite{thurston2} asserts that each mapping class in $\map(S_g)$ is either periodic, reducible, or pseudo-Anosov. Thus, a reducible mapping class $F\in \map(S_g)$ can be decomposed into periodic and pseudo-Anosov components by cutting the surface along the canonical reduction system $\mathcal{C}(F)$ for $F$. This decomposition is known as a \textit{canonical decomposition (or a Nielsen decomposition)}.

A collection of isotopy classes of pairwise disjoint essential (non-nullhomotopic) simple closed curves is called a \textit{multicurve}. Given a multicurve $C$ in $S_g$, the symbol $S_g(C)$ will denote the cut surface (possibly disconnected) obtained by capping the boundary components of $\overline{S_g\setminus N}$ (by marked discs), where $N$ is a closed regular neighborhood of $C$.

A multicurve $C=\{c_1,c_2,\dots,c_m\}$ in $S_g$ is said to be \textit{nonseparating} (or \textit{separating} resp.) if the cut surface $S_g(C)$ is connected (or disconnected resp.). A separating multicurve will be called a \textit{bounding multicurve} if $c_1+c_2+\dots+c_m=0$ in $H_1(S_g;\mathbb{Z})$, that is, the multicurve $C$ bounds a subsurface in $S_g$. For a bounding multicurve $C$, the minimum of genera of components of $S_g(C)$ is said to be the \textit{genus} of the multicurve $C$. A bounding multicurve of size $2$ will be called a \textit{bounding pair}.

An infinite order reducible mapping class is said to be \textit{pseudo-periodic} if there are only periodic components in its canonical decomposition. A nontrivial $G\in \map(S_g)$ is said to be a \textit{root of $F\in \map(S_g)$ of degree $n$} if there exists an integer $n>1$ such that $G^n=F$ and $|G|=n|F|$. A mapping class $F \in \map(S_g)$ is said to be \textit{primitive} if it has no roots of degree $n$ for any $n$.

Given a multicurve $C=\{c_1, c_2,\dots, c_m\}$ in $S_g$ and non-zero integers $q_i$, for $ 1 \leq i \leq m$, a mapping class of the form $T_{c_1}^{q_1}T_{c_2}^{q_2}\cdots T_{c_m}^{q_m}$ is said to be a \textit{multitwist} about $C$. We observe that multitwists are pseudo-periodic mapping classes having trivial periodic components in its canonical decomposition. Therefore, a pseudo-periodic mapping class that is not a multitwist will be called a \textit{nontrivial pseudo-periodic} mapping class. From the Nielsen-Thurston classification, it follows that a nontrivial pseudo-periodic mapping class $F$ is a root of a multitwist about the canonical reduction system $\mathcal{C}(F)$.

\begin{exmp}
\label{exmp:compatibility}
Let $F$ be a root of $T_c$ of degree $n$. Then $F$ is represented by an $\mathcal{F} \in \homeo(S_g)$ such that $\mathcal{F}(N)=N$, where $N$ is a closed annular neighborhood of $c$. Thus, $\mathcal{F}$ induces a $C_n$-action on $S_g(c)$ with two fixed points. Moreover, the sum of induced rotation angles about these fixed points is $2\pi /n$ modulo $2\pi$.  Conversely, given periodic mapping classes having a (two, in case $c$ is nonseparating) distinguished fixed points such that the sum of induced rotation angles about these fixed points is $2\pi /n$ modulo $2\pi$, one can reverse this process to recover the root $F$ of $T_c$. (We refer the reader to~\cite{rajeevsarathy3,rajeevsarathy5} for details.) We illustrate this construction of roots of Dehn twists in Figure \ref{fig2}.
\begin{figure}[ht]
\tiny
\begin{subfigure}{\textwidth}
\centering
\labellist
\pinlabel $(2,5)$ at 202, 85
\pinlabel $(2,5)$ at 202, 35
\pinlabel $D$ at 125, 47
\pinlabel $c$ at 285, 60
\endlabellist
\centering
\includegraphics[scale=0.6]{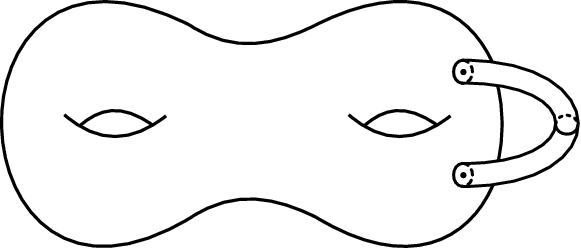}
\caption{The sum of local rotation angles about the two fixed points of cyclic action $D=(5,0;(1,5),(2,5),(2,5))$ associated with the pair $(2,5)$ is $2\pi/5$ modulo $2\pi$. Hence, the action $D$ can be extended to a root of $T_c$ in $\map(S_3)$ of degree $5$, where $c$ is a non-separating curve.}
\label{subfig1}
\end{subfigure}

\begin{subfigure}{\textwidth}
\labellist
\small
\pinlabel $(1,8)$ at 220, 50
\pinlabel $(9,10)$ at 330, 45
\pinlabel $D_1$ at 120, 40
\pinlabel $D_2$ at 425, 40
\pinlabel $c$ at 275, 40
\endlabellist
\centering
\vspace{0.5cm}
\includegraphics[scale=0.6]{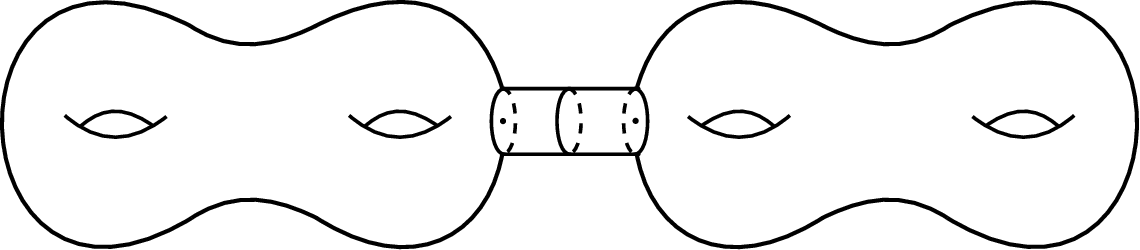}
\caption{The sum of local rotation angles about the two fixed points of cyclic actions $D_1=(8,0;(1,2),(3,8),(1,8))$ and $D_2=(10,0;(1,2),(3,5),(9,10))$ associated with the pairs $(1,8)$ and $(9,10)$, respectively, is $2\pi/40$ modulo $2\pi$. Since $\lcm(8,10)=40$, a root of $T_c$ can be constructed from $D_1$ and $D_2$ in $\map(S_4)$ of degree $40$, where $c$ is a separating curve.}
\label{subfig2}
\end{subfigure}
\caption{Construction of a root of Dehn twist.}
\label{fig2}
\end{figure}
\end{exmp}

The angle sum condition in Example~\ref{exmp:compatibility} (in the construction of the root of a Dehn twist) generalizes to a formal ``compatibility condition" between pairs of orbits of one or more cyclic actions. 

\begin{defn}
\label{defn:twist_factor}
For $i = 1,2$, let $O_i$ be an orbit of cyclic action $D_i$ such that $|O_1| = |O_2|$.
\begin{enumerate}[(i)]
\item We say that $O_1$ and $O_2$ are \textit{trivially compatible} if $|O_1| = |O_2| = n$ (in this case $n(D_1)=n(D_2)$). 
\item Let the pair $(d_i,n_i)$ correspond to the orbit $O_i$ in the data set $D_i$, where we assume that $(d_i,n_i)=(0,1)$ if $|O_i|=n(D_i)$. Let $n=\lcm(n(D_1),n(D_2))$ and \begin{equation}
\label{eqn:twist_compa}
\frac{2\pi d_1^{-1}}{n_1}+\frac{2\pi d_2^{-1}}{n_2}\equiv \frac{2\pi k'}{n} \pmod{2 \pi},
\end{equation}
where $d_id_i^{-1}\equiv 1 \pmod{n_i}$ and $k'$ is some integer. We say that the orbits $O_1$ and $O_2$ are \textit{compatible with twist factor $k$} if $k\equiv k' \pmod n$ and $-n/2<k\leq n/2$.
\end{enumerate}
\end{defn}

\noindent We conclude this section with the following technical remark.
 
\begin{rem}
\label{rem:gcd_remark}
Let $G$ be a root of the multitwist $T_{c_1}^{q_1}T_{c_2}^{q_2}\cdots T_{c_m}^{q_m}$ of degree $n$. Since each $q_i=ns_i+r_i$, where $|r_i|\leq n/2$, it suffices to assume that $-n/2<q_i\leq n/2$ and $q_i\neq 0$ for every $i$. To see this, we write $q_i=ns_i'+r_i'$, where $0\leq r_i'<n$. If $r_i'\leq n/2$, take $s_i=s_i'$ and $r_i=r_i'$. If $r_i'>n/2$, then take $r_i=-(n-r_i')$ and $s_i=s_i'+1$. Hence, for $q_i=ns_i+r_i$, we have
\begin{center}
$G^n=(T_{c_1}^{s_1}T_{c_2}^{s_2}\cdots T_{c_m}^{s_m})^n T_{c_1}^{r_1}T_{c_2}^{r_2}\cdots T_{c_m}^{r_m}$,
\end{center}
where $-n/2<r_i \leq n/2$.
\end{rem}
\noindent From here on we consider only those nontrivial pseudo-periodics which satisfy the condition of Remark \ref{rem:gcd_remark}. 

\section{Primitivity in the mapping class group}
\label{sec:gen_prim}
In this section, we will use the theory of pseudo-periodic mapping classes developed in \cite{rajeevsarathy6, matsumoto} and Thurston's orbifold theory~\cite[Chapter 13]{thurston1} to formulate a combinatorial data set that encodes the conjugacy class of a pseudo-periodic mapping class. Using this encoding and a result from~\cite{dhanwani1} concerning the primitivity of periodics, we derive equivalent conditions for the primitivity of pseudo-periodic mapping classes. We then apply this characterization (and a couple of known algorithms) to obtain an algorithm that determines the primitivity of an arbitrary mapping class up to conjugacy.

\subsection{A combinatorial encoding of pseudo-periodic mapping classes}
\begin{figure}[ht]
\centering
\tiny
\labellist

\pinlabel $D_1$ at 80, 17
\pinlabel $D_2$ at 290, 17
\pinlabel $D_3$ at 500, 17
\pinlabel $c_1$ at 185, 35
\pinlabel $c_2$ at 400, 35

\pinlabel $(1,21)$ at 145, 52
\pinlabel $(19,20)$ at 232, 52
\pinlabel $(11,20)$ at 347, 52
\pinlabel $(9,11)$ at 437, 52

\endlabellist
\includegraphics[scale=0.7]{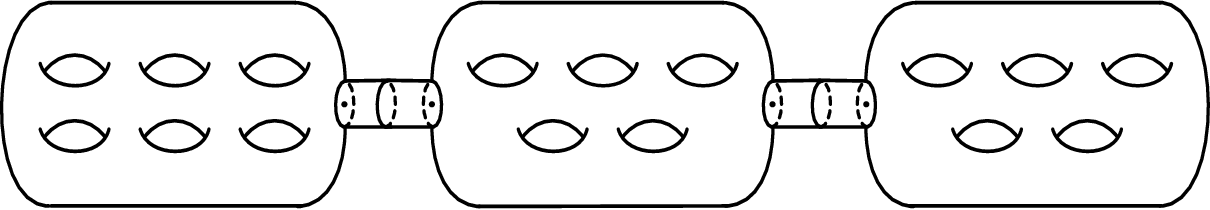}
\caption{An illustration of the tuple $\mathcal{D}$.}
\label{fig3}
\end{figure}

Consider the following cyclic actions
\begin{enumerate}[(i)]
\item $D_1=(21,0;(2,3),(2,7),(1,21)_1)$,
\item $D_2=(20,0;(1,2),(11,20)_2,(19,20)_1)$, and
\item $D_3=(11,0;(1,11),(1,11),(9,11)_2)$
\end{enumerate}
on the surfaces $S_6$, $S_5$, and $S_5$ respectively. We observe that the pair of cone points with the suffix $1$ (or $2$ resp.) are compatible with twist factor $-1$ (or $1$ resp.). The data sets $D_1$, $D_2$, and $D_3$ are periodic components of a pseudo-periodic mapping class $F\in \map(S_{16})$, which is a root of the multitwist $G=T_{c_1}^{-11}T_{c_2}^{21}$ of degree $4620$, where $c_i$ is a separating curve in $S_{16}$ shown in Figure \ref{fig3}. We observe that the canonical reduction system of $F$ is $\mathcal{C}(F)=\{c_1,c_2\}$.

\noindent The information associated with the conjugacy class of the pseudo-periodic mapping class $F$ can be encoded as
\begin{center}
$\mathcal{D}=\llbracket 4620;(D_1,1),(D_2,1),(D_3,1);\llparenthesis 1,-11,-1;1,2\rrparenthesis,\llparenthesis 1,21,1;2,3\rrparenthesis;-\rrbracket$.
\end{center}
In the tuple $\mathcal{D}$, the integer $4620$ is the degree of the root and the pairs $(D_i,1)$ encode the fact that each $D_i$ acts on a surface orbit of size $1$ induced by $F$ on the cut surface $S_{16}(\mathcal{C}(F))$. We observe that the orbit of the curve $c_1$ (resp. $c_2$) is of size $1$ under the action of $F$ on $\mathcal{C}(F)$ with twist factor $-1$ (resp. $1$). An annular neighborhood of $c_1$ (resp. $c_2$) connects components of $S_{16}(\mathcal{C}(F))$ on which $D_1$ and $D_2$ (resp. $D_2$ and $D_3$) act, and the exponent of $T_{c_1}$ (resp. $T_{c_2}$) in the multitwist $G$ is $-11$ (resp. $21$). The tuple $\llparenthesis 1,-11,-1;1,2\rrparenthesis$ (resp. $\llparenthesis 1,21,1;2,3\rrparenthesis$) encodes the information associated with the curve $c_1$ (resp. $c_2$) as we have just observed. This motivates the following definition.

\begin{defn}
\label{defn:encoding}
For $1\leq i,i_r,j_r \leq s$, $1\leq r\leq \ell$, and $\ell_1\leq \ell_2\leq \ell_3\leq \ell$, a \textit{pseudo-periodic data set} is a tuple of the form
\begin{align*}
\mathcal{D}=&\llbracket n;(D_1,p_1),\dots,(D_s, p_s);\llparenthesis m_1,q_1,k_1;i_1,j_1\rrparenthesis,\dots,\llparenthesis m_{\ell_1},q_{\ell_1},k_{\ell_1};i_{\ell_1},j_{\ell_1}\rrparenthesis,\llparenthesis m_{\ell_1+1},q_{\ell_1+1},k_{\ell_1+1};\\&i_{\ell_1+1},j_{\ell_1+1}\rrparenthesis^*,\dots,\llparenthesis m_{\ell_2},q_{\ell_2},k_{\ell_2};i_{\ell_2},j_{\ell_2}\rrparenthesis^*;(m_{\ell_2+1},q_{\ell_2+1},k_{\ell_2+1};i_{\ell_2+1},j_{\ell_2+1}),\dots,(m_{\ell_3},q_{\ell_3},k_{\ell_3};\\&i_{\ell_3},j_{\ell_3}),(m_{\ell_3+1},q_{\ell_3+1},k_{\ell_3+1};i_{\ell_3+1},j_{\ell_3+1})^*,\dots,(m_{\ell},q_{\ell},k_{\ell};i_{\ell},j_{\ell})^*\rrbracket,
\end{align*}
where $n$, $p_i$, and $m_r$ are positive integers, $q_r$ and $k_r$ are non-zero integers satisfying the following conditions.
\begin{enumerate}[(i)]
\item For $1\leq i\leq s$, each $D_i$ is a cyclic data set of degree $n(D_i)$.
\item For $1\leq r\leq \ell$, $i_r\leq j_r$, $p_{i_r}\leq p_{j_r}$, $\sgn k_r=\sgn q_r$, $k_r \mid q_r$.
\item For $n=1$, $n(D_i)=p_i=m_r=1$, $q_r=k_r$ for all $i,r$, and $\ell_1=\ell_2$, $\ell_3=\ell$. 
\item When $n>1$ following conditions hold.
\begin{enumerate}[(a)] 
\item For $\ell_2+1\leq r\leq \ell_3$, $i_r=j_r$.
\item For each $1\leq r\leq \ell_1$ or $\ell_2+1\leq r\leq \ell_3$, following holds.
\begin{enumerate}[(1)]
\item If $p_{j_r}=1$, then either $n(D_{i_r})=m_r=n(D_{j_r})$ or each of $D_{i_r}$ and $D_{j_r}$ have compatible orbits of size $m_r$ with twist factor $k_r$.
\item If $1=p_{i_r}<p_{j_r}$, then $m_r=p_{j_r}$ and $D_{i_r}^{p_{j_r}}$ has $m_r$ fixed points, corresponding to an orbit of size $m_r$ of $D_{i_r}$, each of which is compatible to a fixed point of $D_{j_r}$ with twist factor $k_r$.
\item If $p_{i_r},p_{j_r}>1$, then $p_{i_r}=m_r=p_{j_r}$, and each of $D_{i_r}$ and $D_{j_r}$ have a fixed point which are compatible with twist factor $k_r$. 
\end{enumerate}
\item For each $\ell_3+1\leq r\leq \ell$, $p_{i_r}=1$, $m_r=1$, $i_r=j_r$, and $D_{i_r}^2$ has $2$ fixed points, corresponding to an orbit of size $2$ of $D_{i_r}$, which are compatible with twist factor $k_r$.
\item For each $\ell_1+1\leq r\leq \ell_2$, $p_{i_r}=2$, $m_r=1$, $i_r=j_r$, and $D_{i_r}$ has a fixed point which is compatible to itself with twist factor $k_r$.
\item For each $1\leq r\leq \ell_1$ and $\ell_2<r\leq \ell_3$, $-n_r/2<k_r\leq n_r/2$, $-n/2<q_r\leq n/2$, and $n=(q_r/k_r)p_{j_r}n_r$, where $n_r=p_{i_r}\lcm (n(D_{i_r})/p_{j_r},n(D_{j_r})/p_{i_r})$.
\item For each $\ell_1< r\leq \ell_2$ and $\ell_3<r\leq \ell$, $-n(D_{i_r})/2<k_r\leq n(D_{i_r})/2$, $-n/2<q_r\leq n/2$, and $n=(q_r/k_r)n(D_{i_r})$.
\end{enumerate}
\end{enumerate}
The positive integer $\min\{|q_1|,|q_2|,\dots,|q_{\ell}|\}$ associated to the pseudo-periodic data set $\mathcal{D}$ will be denoted by $q(\mathcal{D})$.
\end{defn}

\noindent We further explain the notation in Definition~\ref{defn:encoding} with the help of the following examples.

\begin{exmp}
For $i=1,2,3,4$, we construct a pseudo-periodic mapping classes $F_i\in \map(S_4)$ whose canonical reduction systems differ from each other topologically. Let $F$ be the free rotation of $S_4$ by $2\pi /3$ as shown in Figure~\ref{fig4}.
\begin{figure}[ht]
\begin{subfigure}[b]{0.46\textwidth}

\labellist
\small
\pinlabel $c_1$ at 70, 20
\pinlabel $c_2$ at 400, 20
\pinlabel $c_3$ at 260, 345
\pinlabel $D_1$ at 180, 120
\pinlabel $\frac{2\pi}{3}$ at 60, 255
\endlabellist
\centering
\includegraphics[scale=0.43]{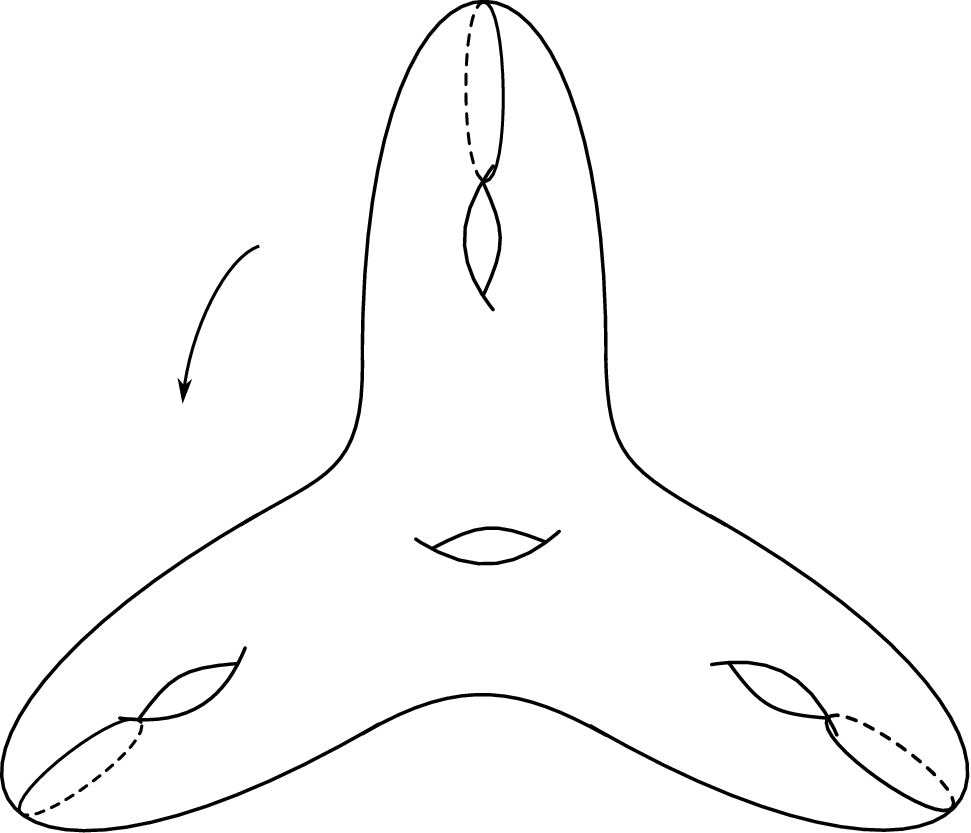}
\vspace{2mm}
\caption{The mapping class $F_1=F\circ T_{c_1}$ is a pseudo-periodic mapping class such that $F_1^3=T_{c_1}T_{c_2}T_{c_3}$. The pseudo-periodic data set corresponding to $F_1$ can be written as $\mathcal{D}_1=\llbracket 3;(D_1,1);-;(3,1,1;1,1)\rrbracket$, where $D_1=(3,1;-)$ is a free $C_3$-action on $S_1$. We observe that $\mathcal{C}(F_1)$ is a nonseparating multicurve.}
\label{subfig:free_action}
\end{subfigure}
\hfill
\begin{subfigure}[b]{0.46\textwidth}
\labellist
\small
\pinlabel $D_1$ at 45, 25
\pinlabel $D_1$ at 425, 20
\pinlabel $D_1$ at 240, 365
\pinlabel $c_1$ at 180, 145
\pinlabel $c_2$ at 265, 100
\pinlabel $c_3$ at 250, 190
\pinlabel $\frac{2\pi}{3}$ at 60, 255
\endlabellist
\centering
\includegraphics[scale=.4]{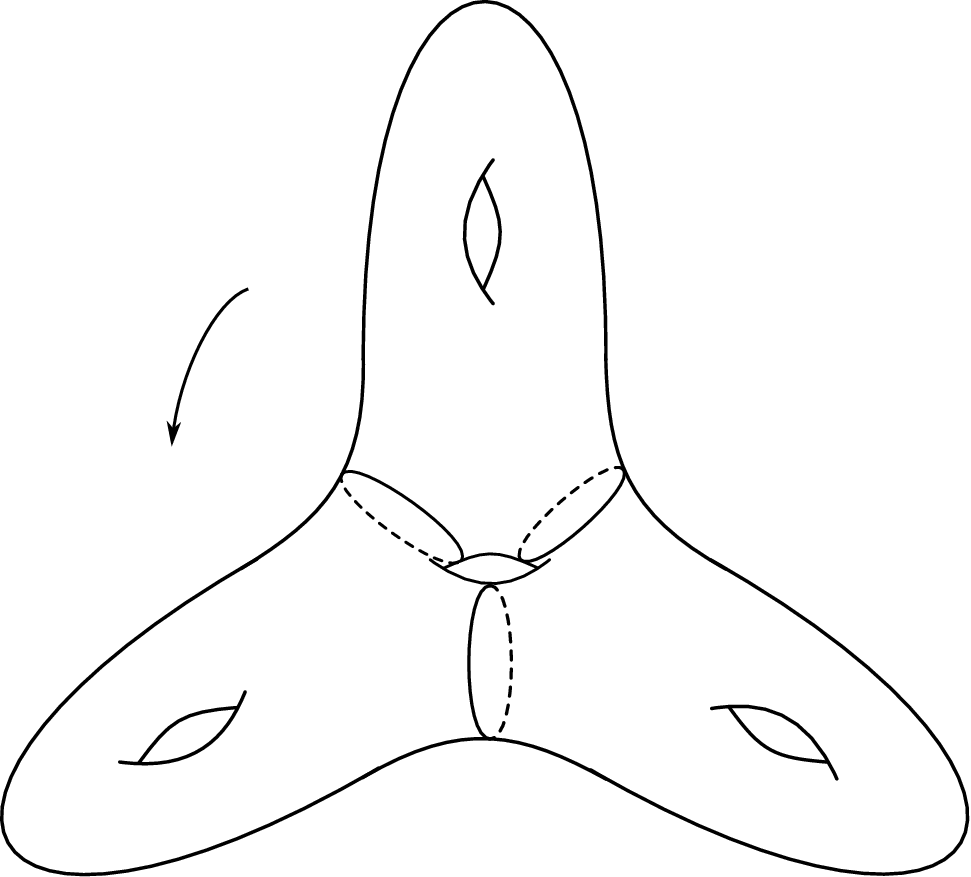}
\vspace{3mm}
\caption{The mapping class $F_2=F\circ T_{c_1}$ is a pseudo-periodic mapping class such that $F_2^3=T_{c_1}T_{c_2}T_{c_3}$. The pseudo-periodic data set corresponding to $F_2$ can be written as $\mathcal{D}_2=\llbracket 3;(D_1,3);\llparenthesis 3,1,1;1,1 \rrparenthesis;-\rrbracket$, where $D_1$ is trivial on $S_1$. We observe that $\mathcal{C}(F_2)$ is a separating multicurve containing only nonseparating curves.}
\label{subfig:free_action2}
\end{subfigure}
\hfill
\begin{subfigure}[b]{0.46\textwidth}
\labellist
\small
\pinlabel $D_2$ at 45, 25
\pinlabel $D_2$ at 425, 20
\pinlabel $D_2$ at 240, 365
\pinlabel $D_1$ at 190, 160
\pinlabel $c_1$ at 180, 120
\pinlabel $c_2$ at 275, 120
\pinlabel $c_3$ at 240, 190
\pinlabel $\frac{2\pi}{3}$ at 60, 255
\endlabellist
\centering
\vspace{7mm}
\includegraphics[scale=0.4]{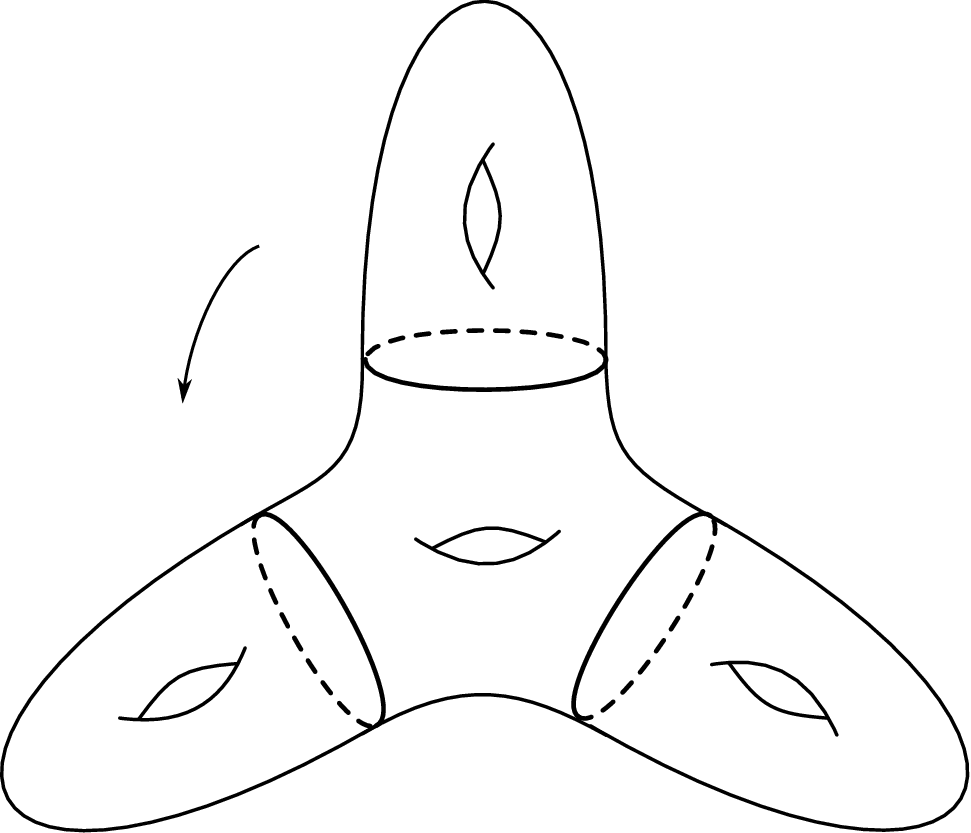}
\vspace{2mm}
\caption{The mapping class $F_3=F\circ T_{c_1}$ is a pseudo-periodic mapping class such that $F_3^3=T_{c_1}T_{c_2}T_{c_3}$. The pseudo-periodic data set corresponding to $F_3$ can be written as $\mathcal{D}_3=\llbracket 3;(D_1,1),(D_2,3);\llparenthesis 3,1,1;1,2 \rrparenthesis;-\rrbracket$, where $D_1=(3,1;-)$ is a free $C_3$-action and $D_2$ is trivial, both on $S_1$. We observe that $\mathcal{C}(F_3)$ is a separating multicurve containing only separating curves of genus $1$.}
\label{subfig:free_action3}
\end{subfigure}
\hfill
\begin{subfigure}[b]{0.46\textwidth}
\labellist
\small
\pinlabel $D_1$ at 40, 10
\pinlabel $D_2$ at 40, 101
\pinlabel $c_1$ at 20, 52
\pinlabel $c_2$ at 82, 52
\pinlabel $c_3$ at 43, 75
\pinlabel $\frac{2\pi}{3}$ at 65, 135
\endlabellist
\centering
\includegraphics[scale=1.3]{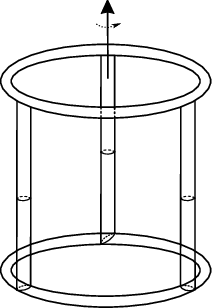}
\vspace{2mm}
\caption{The mapping class $F_4=F\circ T_{c_1}$ is a pseudo-periodic mapping class such that $F_4^3=T_{c_1}T_{c_2}T_{c_3}$. The pseudo-periodic data set corresponding to $F_4$ can be written as $\mathcal{D}_4=\llbracket 3;(D_1,1),(D_2,1);\llparenthesis 3,1,1;1,2 \rrparenthesis;-\rrbracket$, where for $i=1,2$, $D_i=(3,1;-)$ is a free $C_3$-action on $S_1$. We observe that $\mathcal{C}(F_4)$ is a bounding multicurve of genus $1$.}
\label{subfig:free_action4}
\end{subfigure}
\caption{Pseudo-periodic mapping classes in $\map(S_4)$ with distinct canonical reduction system.}
\label{fig4}
\end{figure}
\end{exmp}

\begin{exmp}
In this example, we describe a side-exchanging root of a multitwist in $\map(S_3)$ and $\map(S_4)$ as shown in Figure~\ref{fig5}. In the pseudo-periodic data sets (encoding the roots), we put a $*$ in superscript of a tuple to indicate that it corresponds to an orbit of curves whose sides get exchanged.

\begin{figure}[ht]
\begin{subfigure}{0.46\textwidth}
\labellist
\small
\pinlabel $(1,3)$ at 200, 85
\pinlabel $(1,3)$ at 200, 35
\pinlabel $D_1$ at 125, 47
\pinlabel $c$ at 285, 60
\endlabellist
\centering
\includegraphics[scale=0.7]{selfpp}
\caption{Consider the cyclic action $D_1=(6,0;(1,3),(5,6),(5,6))$ on $S_2$. We observe that the cone point of order $3$ corresponds to a distinguished orbit of size $2$ on $S_2$. These points are compatible with twist factor $-2$. Therefore, a side-exchanging root $F$ of $T_c^{-2}$ of degree $6$ can be constructed from $D_1$, where $c$ is a nonseparating curve. The pseudo-periodic data set corresponding to this root can be written as $\mathcal{D}=\llbracket 6;(D_1,1);-;(1,-2,-2;1,1)^*\rrbracket$.}
\end{subfigure}
\hfill
\begin{subfigure}{0.52\textwidth}
\labellist
\small
\pinlabel $c_3$ at 105, 10
\pinlabel $c_2$ at 105, 45
\pinlabel $c_1$ at 105, 80
\pinlabel $D_1$ at 50, 20
\pinlabel $D_1$ at 180, 20
\endlabellist
\centering
\includegraphics[scale=.9]{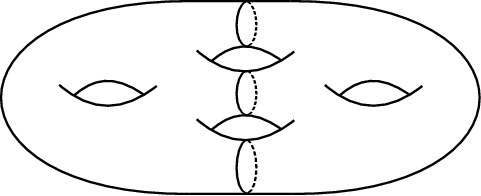}
\caption{Consider the cyclic action $D_1=(3,0;(2,3),(2,3),(2,3))$ on $S_1$. We take two copies of $S_1$ with cyclic action $D_1$. We observe that any two fixed points of $D_1$ are compatible with twist factor $1$. We can construct a side-exchanging root $F$ of $T_{c_1}^2T_{c_2}^2T_{c_3}^2$ of degree $6$ with periodic component $D_1$. The pseudo-periodic data set corresponding to $F$ can be written as $\mathcal{D}=\llbracket 6;(D_1,2);\llparenthesis 1,2,1;1,1 \rrparenthesis^*,\llparenthesis 1,2,1;1,1 \rrparenthesis^*,\llparenthesis 1,2,1;1,1 \rrparenthesis^*;-\rrbracket$.}
\end{subfigure}
\caption{Side-exchanging roots of multitwists in $\map(S_3)$ and $\map(S_4)$.}
\label{fig5}
\end{figure}
\end{exmp}

\noindent The following definition states permissible reordering of tuples allowed in a pseudo-periodic data set.
\begin{defn}
The pseudo-periodic data sets
\begin{align*}
\mathcal{D}=&\llbracket n;(D_1,p_1),\dots,(D_s, p_s);\llparenthesis m_1,q_1,k_1;i_1,j_1\rrparenthesis,\dots,\llparenthesis m_{\ell_1},q_{\ell_1},k_{\ell_1};i_{\ell_1},j_{\ell_1}\rrparenthesis,\llparenthesis m_{\ell_1+1},q_{\ell_1+1},k_{\ell_1+1};\\&i_{\ell_1+1},j_{\ell_1+1}\rrparenthesis^*,\dots,\llparenthesis m_{\ell_2},q_{\ell_2},k_{\ell_2};i_{\ell_2},j_{\ell_2}\rrparenthesis^*;(m_{\ell_2+1},q_{\ell_2+1},k_{\ell_2+1};i_{\ell_2+1},j_{\ell_2+1}),\dots,(m_{\ell_3},q_{\ell_3},k_{\ell_3};\\&i_{\ell_3},j_{\ell_3}),(m_{\ell_3+1},q_{\ell_3+1},k_{\ell_3+1};i_{\ell_3+1},j_{\ell_3+1})^*,\dots,(m_{\ell},q_{\ell},k_{\ell};i_{\ell},j_{\ell})^*\rrbracket,
\end{align*}
and
\begin{align*}
\mathcal{D}'=&\llbracket n';(D_1',p_1'),\dots,(D_{s'}', p_{s'}');\llparenthesis m_1',q_1',k_1';i_1',j_1'\rrparenthesis,\dots,\llparenthesis m_{\ell_1'}',q_{\ell_1'}',k_{\ell_1'}';i_{\ell_1'}', j_{\ell_1'}'\rrparenthesis,\llparenthesis m_{\ell_1'+1}',q_{\ell_1'+1}',k_{\ell_1'+1}';\\&i_{\ell_1'+1}',j_{\ell_1'+1}'\rrparenthesis^*,\dots,\llparenthesis m_{\ell_2'}',q_{\ell_2'}',k_{\ell_2'}';i_{\ell_2'}',j_{\ell_2'}'\rrparenthesis^*;(m_{\ell_2'+1}',q_{\ell_2'+1}',k_{\ell_2'+1}';  i_{\ell_2'+1}',j_{\ell_2'+1}'), \dots,(m_{\ell_3'}',q_{\ell_3'}',k_{\ell_3'}';\\&i_{\ell_3'}', j_{\ell_3'}'),(m_{\ell_3'+1}',q_{\ell_3'+1}',k_{\ell_3'+1}';i_{\ell_3'+1}',j_{\ell_3'+1}')^*,\dots,(m_{\ell'}',q_{\ell'}',k_{\ell'}';i_{\ell'}',j_{\ell'}')^*\rrbracket,
\end{align*}
are said to be \textit{equivalent} is the following conditions hold.
\begin{enumerate}[(i)]
\item $n=n'$, $s=s'$, $\ell=\ell'$, and $\ell_i=\ell_i'$ for $i=1,2,3$.
\item $\{(D_1,p_1),\dots,(D_s,p_s)\}=\{(D_1',p_1'),\dots,(D_{s'}',p_{s'}')\}$.
\item $\!
\begin{aligned}[t]
\{\llparenthesis m_1,q_1,k_1;i_1,j_1 \rrparenthesis,\dots ,\llparenthesis m_{\ell_1},q_{\ell_1},k_{\ell_1};i_{\ell_1},j_{\ell_1}& \rrparenthesis\}  = \\
& \{\llparenthesis m_1',q_1',k_1';i_1',j_1'\rrparenthesis,\dots, \llparenthesis m_{\ell_1'}',q_{\ell_1'}',k_{\ell_1'}';i_{\ell_1'}', j_{\ell_1'}'\rrparenthesis\}.
\end{aligned}$
\item $\!
\begin{aligned}[t]
\{\llparenthesis m_{\ell_1+1},q_{\ell_1+1},k_{\ell_1+1};i_{\ell_1+1},&j_{\ell_1+1}\rrparenthesis^*,\dots,\llparenthesis m_{\ell_2},q_{\ell_2},k_{\ell_2};i_{\ell_2},j_{\ell_2}\rrparenthesis^*\}=\\
&\{\llparenthesis m_{\ell_1'+1}',q_{\ell_1'+1}',k_{\ell_1'+1}';i_{\ell_1'+1}',j_{\ell_1'+1}'\rrparenthesis^*,\dots,\llparenthesis m_{\ell_2'}',q_{\ell_2'}',k_{\ell_2'}';i_{\ell_2'}',j_{\ell_2'}'\rrparenthesis^*\}.
\end{aligned}$
\item $\!
\begin{aligned}[t]
\{(m_{\ell_2+1},
q_{\ell_2+1},k_{\ell_2+1};&i_{\ell_2+1},j_{\ell_2+1}),\dots,(m_{\ell_3},q_{\ell_3},k_{\ell_3};i_{\ell_3},j_{\ell_3})\}=\\
&\{(m_{\ell_2'+1}',q_{\ell_2'+1}',k_{\ell_2'+1}';i_{\ell_2'+1}', j_{\ell_2'+1}'),\dots,(m_{\ell_3'}',q_{\ell_3'}',k_{\ell_3'}';i_{\ell_3'}', j_{\ell_3'}')\}.
\end{aligned}$
\item $\!
\begin{aligned}[t]
\{(m_{\ell_3+1},q_{\ell_3+1},k_{\ell_3+1};i_{\ell_3+1},&j_{\ell_3+1})^*,\dots,(m_{\ell},q_{\ell},k_{\ell};i_{\ell},j_{\ell})^*\}=\\
&\{(m_{\ell_3'+1}',q_{\ell_3'+1}',k_{\ell_3'+1}';i_{\ell_3'+1}',j_{\ell_3'+1}')^*,\dots,(m_{\ell'}',q_{\ell'}',k_{\ell'}';i_{\ell'}',j_{\ell'}')^*\}.
\end{aligned}$
\end{enumerate}
\end{defn}

\noindent The following basic property of multitwists plays crucial role in our theory. 
\begin{lem}[{\cite[Lemma 3.17]{primer}}]
\label{lem:equal_multitwist}
Let $A=\{a_1,a_2,\dots,a_n\}$ and $B=\{b_1,b_2,\dots,b_m\}$ be two multicurves in a surface $S$. Let $p_i$ and $q_i$ be non-zero integers. If
\begin{center}
$T_{a_1}^{p_1}T_{a_2}^{p_2}\cdots T_{a_m}^{p_n}=T_{b_1}^{q_1}T_{b_2}^{q_2}\cdots T_{b_m}^{q_m}$
\end{center}
in $\map(S)$, then $m=n$ and the sets $\{T_{a_i}^{p_i}\}$ and $\{T_{b_i}^{q_i}\}$ are equal.
\end{lem}

\noindent From Lemma~\ref{lem:equal_multitwist}, it follows that a root of a multitwist about multicurve $C$ preserves $C$.

\begin{lem}
\label{lem:root_preserve}
Let $G$ be a root of a multitwist $F=T_{c_1}^{q_1}T_{c_2}^{q_2}\cdots T_{c_m}^{q_m}$ about multicurve $C=\{c_1,c_2,\dots,c_m\}$. Then, $G$ can be modified by an isotopy so that it preserves $C$.
\begin{proof}
Since $G^n=F$, we have $F=GFG^{-1}$ and therefore,
\begin{center}
$T_{c_1}^{q_1}T_{c_2}^{q_2}\cdots T_{c_m}^{q_m}=T_{G(c_1)}^{q_1}T_{G(c_2)}^{q_2}\cdots T_{G(c_m)}^{q_m}$. 
\end{center}
It follows from Lemma \ref{lem:equal_multitwist} that, $G(C)=C$.
\end{proof}
\end{lem}

\begin{rem}
\label{rem:cyclic_multiple}
Let $S_g(p)$ denote the disjoint union of $p$-copies of surface $S_g$, that is, $S_g(p)=\{1,2,\dots,p\}\times S_g$. The component $\{i\}\times S_g$ will be denoted by $S_g^i$. Let $F\in \map(S_g)$ be a mapping class extended by the identity on $\{1,2,\dots,p\}$. Let $\sigma_p$ be a $p$-cycle on $\{1,2,\dots,p\}$ extended by the identity on $S_g$. Then $\sigma_pF\in \map(S_g(p))\cong \oplus_{i=1}^p \map(S_g^i)$. Since
\[
(\sigma_pF)^p=(\sigma_p F\sigma_p^{-1})(\sigma_p^2 F\sigma_p^{-2})\cdots (\sigma_p^{p-1} F\sigma_p^{1-p})F,
\]
each $\sigma_p^{i}F \sigma_p^{-i}\in \map(S_g^i)$ called the \textit{first return map} on $S_g^i$. Note that all of these first return maps are pairwise conjugate to each other. If $F$ is periodic of order $n$, then $\sigma_pF$ is also periodic of order $np$. Conversely, given $\tilde{F}\in \map(S_g(p))$, for $F:=(\tilde{F})^p\big|_{S_g^1}$, we have $F\in \map(S_g^1)$ and $\tilde{F}=\sigma_p F$. Furthermore, $(\tilde{F})^p\big|_{S_g^i}\in \map(S_g^i)$ are all pairwise conjugate to each other. If $\tilde{F}$ is periodic of order $n$, then $F$ is periodic of order $n/p$. Hence, there is a one-to-one correspondence between $C_n$-actions on $S_g(p)$ and $C_{n/p}$-actions on $S_g$.
\end{rem}

\noindent The following proposition is a complete generalization of the theory developed in \cite{rajeevsarathy6}.
\begin{prop}
\label{thm:main_theorem}
There is a correspondence between equivalence classes of pseudo-periodic data sets and the conjugacy classes of pseudo-periodic mapping classes.
\begin{proof}We will follow the notations of the Definition \ref{defn:encoding}. Let $G$ be a pseudo-periodic mapping class. Assume that $G$ is a root of a multitwist of degree $n$, that is, $G^n=T_{C_1}^{q_1}T_{C_2}^{q_2}\cdots T_{C_{\ell}}^{q_{\ell}}$ (take $n=1$ in the case $G$ is itself a multitwist), where $C_r$ is an orbit of curves of size $m_r$ under the action of $G$ on $C=\sqcup_{r=1}^{\ell}C_r$ and $T_{C_r}$ is the product of Dehn twists about the curves in the orbit $C_r$ (from Lemma \ref{lem:root_preserve} $G$ preserves multicurve $C$). Let $N$ be a closed regular neighborhood of $C$. We can change $G$ by isotopy so that $G(N)=N$. Then $G$ acts on the components of $N$ and $S_g(C)$, and induces a $C_n$-action on $S_g(C)$.

Assume that $S_g(C)=\sqcup_{i=1}^s S_{g_i}(p_i)$, where $S_{g_i}(p_i)$ is a \textit{surface orbit} under the action of $G$ on $S_g(C)$. Then $G$ induces a $C_{u_i}$-action on $S_{g_i}(p_i)$ for some positive integer $u_i$ such that $u_i\mid n$ for each $i$. Let $D_i$ be the cyclic data set corresponding to $C_{u_i/p_i}$-action (that is, the first return map) induced by $C_{u_i}$-action on $S_{g_i}$ (see Remark \ref{rem:cyclic_multiple}). Assume that $N=\sqcup_{r=1}^{\ell} N_r$, where $N_r$ is an orbit of annuli (that is, regular neighborhood of curves in $C_r$) under the action of $G$ on $N$. Let the orbit $N_r$ connects the surface orbits corresponding to the indices $i_r$ and $j_r$, where $1\leq i_r,j_r \leq s$. Without loss of generality, we assume that $p_{i_r}\leq p_{j_r}$.

When $G$ is a multitwist, all periodic components are trivial, and so $n(D_i)=1$ for every $i$. Since the multitwist $G$ fixes each component of $S_g(C)$ and $C$, we have $p_i=m_r=1$ and $k_r=q_r$ for every $i,r$. Now, assume that $G$ is not a multitwist. The action of $G^{p_{j_r}}$ twists the annuli in the orbit $N_r$ by $2\pi k_r/n_r$, for some non-zero integer $k_r$ such that $-n_r/2<k_r\leq n_r/2$, where
\begin{center}
$n_r=p_{i_r}\lcm \left(\frac{n(D_{i_r})}{p_{j_r}},\frac{n(D_{j_r})}{p_{i_r}}\right)$.
\end{center}
Since $G^n$ twists the annuli in the orbit $N_r$ by $2\pi q_r$, the action of $G^{p_{j_r}}$ must twist these annuli by $(2\pi q_r p_{j_r})/n$. Therefore, we have $q_r/k_r=n/(p_{j_r}n_r)$, where $p_{j_r}n_r=\lcm(u_{i_r},u_{j_r})$. It follows that $k_r\mid q_r$ and $p_{j_r}n_r=\frac{n}{q_r/k_r}$.

First we assume that $G$ preserves the sides of annuli in the orbit $N_r$. If $p_{j_r}=1$, then $D_{i_r}$ and $D_{j_r}$ have distinguished orbits $O_{i_r}$ and $O_{j_r}$, respectively, of size $m_r$ induced by $N_r$. If $m_r=\lcm(n(D_{i_r}),n(D_{j_r}))$, then $n(D_{i_r})=m_r=n(D_{j_r})$. Otherwise, the orbits $O_{i_r}$ and $O_{j_r}$ are compatible with twist factor $k_r$. If $p_{i_r}=1$ and $p_{j_r}>1$, then we must have $p_{j_r}=m_r$, and $D_{i_r}^{p_{j_r}}$ has $m_r$ fixed points, corresponding to an orbit of size $m_r$ of $D_{i_r}$, each of which is compatible to with a fixed point of $D_{j_r}$ with twist factor $k_r$. If both $p_{i_r},p_{j_r}>1$, then we must have $p_{i_r}=p_{j_r}=m_r$, and each of $D_{i_r}$ and $D_{j_r}$ have a fixed point which are compatible with twist factor $k_r$.

Now, we assume that $G$ exchanges the sides of annuli in the orbit $N_r$. We must have that $m_r=1$, $i_r=j_r$, and $p_{i_r}=2$ (or $1$) depending upon whether $C_r$ contained in a submulticurve $C'$ of $C$ which bounds a subsurface in $S_g$ (or not). When $p_{i_r}=1$, since $G$ interchanges the two sides of $N_r$, $D_{i_r}$ has an orbit of size $2$ such that the corresponding fixed points of $D_{i_r}^2$ are compatible with twist factor $k_r$. When $p_{i_r}=2$, since there are two copies of $S_{g(D_{i_r})}$, $N_r$ induces a fixed point on each copy of $S_{g(D_{i_r})}$ both of which are compatible with twist factor $k_r$.

If $C'$ is a bounding multicurve whose sides are exchanged by $G$, then there exists an involution $F\in \map(S_g)$ such that $GF$ is side-preserving. Since $(GF)^2=G(FGF^{-1})$, the degree of $GF$ is $2n$ and the exponents of Dehn twists appearing in the multitwist $(GF)^{2n}$ are $2q_1,2q_2\dots 2q_{\ell}$. Now it follows that $2n=(2q_r/k_r)n(D_{i_r})$, that is, $n=(q_r/k_r)n(D_{i_r})$.

A pseudo-periodic data set $\mathcal{D}_G$ will be constructed inductively as follows. A union of orbits of curves that disconnect the surface (but no strictly smaller union of orbits disconnect the surface) will be encoded by the tuple $\llparenthesis m_r,q_r,k_r;i_r,j_r \rrparenthesis $ (resp. $\llparenthesis 1,q_r,k_r;i_r,j_r \rrparenthesis^*$) depending upon whether $G$ preserves (resp. does not preserve) the sides of curves in the union, and we assume that the total number of such orbits are $\ell_1$ (resp. $\ell_2-\ell_1$). The remaining orbits of curves will be encoded by the tuple $(m_r,q_r,k_r;i_r,j_r)$ (resp. $(1,q_r,k_r;i_r,j_r)^*$) depending upon whether $G$ preserves (resp. does not preserve) the sides of the curves, and we assume that the total number of such orbits are $\ell_3-\ell_2$ (resp. $\ell-\ell_3$). In the case $n=1$, that is, $G$ is a multitwist, $G$ preserves sides of all curves, and therefore $\ell_1=\ell_2, \ell_3=\ell$. The cyclic actions on the surface orbits will be encoded by the pair $(D_i,p_i)$. Hence, we obtain a pseudo-periodic data set $\mathcal{D}_G$ associated with the mapping class $G$.

If $G'=HGH^{-1}$, then $G$ and $G'$ induce the same orbit structure on $C$ and $S_g(C)$, and $H$ maps the orbits induced by $G$ to the orbits induced by $G'$. Thus, $G$ and $G'$ induce equivalent pseudo-periodic data sets.

Conversely, we assume that a pseudo-periodic data set $\mathcal{D}$ is given. Let $G_i\in \map(S_{g(D_i)})$ be a periodic mapping class corresponding to the data set $D_i$ of order $n(D_i)$. Consider the surface $\sqcup_{i=1}^s S_{g(D_i)}(p_i)$. For $1\leq r\leq \ell$, let $N_r$ (with associated twist factor $k_r$) be a disjoint union of annuli of size $m_r$. First we assume that $r\in \{1,2,\dots, \ell_1, \ell_2+1, \dots, \ell_3\}$. If $p_{j_r}=1$, then we know that each of $D_{i_r}$ and $D_{j_r}$ has a distinguished orbit $O_{i_r}$ and $O_{j_r}$ of size $m_r$, respectively, which are either trivially compatible or compatible with twist factor $k_r$. Remove disjoint open disks around points in the orbits $O_{i_r}$ and $O_{j_r}$, and attach the annuli from $N_r$ between $S_{g(D_{i_r})}$ and $S_{g(D_{j_r})}$ with full $k_r$-twists (resp. with $2\pi k_r/n_r$-twists) depending upon whether the orbits were trivially compatible (resp. compatible with twist factor $k_r$). It must be noted that a positive (resp. negative) value of twist factor corresponds to a left-handed (resp. right-handed) twisting of the annuli. We perform a similar construction around compatible fixed points when $p_{j_r}>1$, and when $r\in \{\ell_1+1,\dots,\ell_2,\ell_3+1,\dots,\ell\}$. After performing this construction for all $r\in \{1,2,\dots,\ell\}$, we obtain a pseudo-periodic mapping class $G$ such that $\mathcal{D}_G=\mathcal{D}$.

Now, let $G$ and $G'$ be two mapping classes constructed (as described above) with equivalent pseudo-periodic data sets $\mathcal{D}_G$ and $\mathcal{D}_{G'}$, respectively. Then the orbit structure of $\mathcal{D}_G$ and $\mathcal{D}_{G'}$ are same, up to reordering. Therefore, without loss of generality, let $D_i=D_i'$, where the $D_i$ and $D_i'$ are periodic components of $\mathcal{D}_G$ and $\mathcal{D}_{G'}$, respectively and $1\leq i\leq s$. This implies that periodic mapping classes corresponding to $D_i$ and $D_i'$ are conjugate by some $H_i\in \map(S_{g(D_i)})$ for $1\leq i\leq s$. A reducible mapping class $H\in \map(S_g)$ can be easily constructed having periodic components $H_i$ such that $G'=HGH^{-1}$. Hence, up to conjugacy, $G^n=T_{C_1}^{q_1}T_{C_2}^{q_2}\cdots T_{C_{\ell}}^{q_{\ell}}$, where $C_r$ is a multicurve contained in the orbit $N_r$ for $1\leq r \leq \ell$.
\end{proof}
\end{prop}

From here on, the conjugacy class of a pseudo-periodic mapping class $F$ will be represented by its corresponding pseudo-periodic data set $\mathcal{D}_F$. From Proposition \ref{thm:main_theorem}, it follows that there exists a unique positive integer $q(F)=q(\mathcal{D}_F)$ associated to the pseudo-periodic $F$.

\subsection{Primitivity of pseudo-periodic mapping classes}
Given two conjugacy classes $\mathcal{D}$ and $\mathcal{D}'$ of pseudo-periodic mapping classes, we say $\mathcal{D}$ \textit{carries a root of $\mathcal{D}'$ of degree $n$} if there exists representatives $F$ and $G$ of $\mathcal{D}$ and $\mathcal{D}'$ respectively such that $F^n = G$ in $\map(S_g)$. We will need the following result from~\cite[Proposition 5.1]{dhanwani1} that characterizes the primitivity of periodic mapping classes up to conjugacy.

\begin{prop}
\label{prop:periodic_primitive}
For $g \geq 2$, let $F\in \map(S_g)$ be a periodic mapping class with
\begin{center}
 $D_F=(n,g_0,r;(c_1,n_1),(c_2,n_2),\dots,(c_{\ell},n_{\ell}))$.
\end{center}
Then $F$ has a root $G$ of degree $m$ if and only if the following conditions hold.
\begin{enumerate}[(i)]
\item There exists a homeomorphism $\mathcal{G'}$ of $S_{g_0}$ with
\begin{center}
$D_{G'} = (m,g',r';(d_1,m_1),(d_2,m_2),\dots,(d_k,m_k))$
\end{center}
which induces an automorphism (preserve the size of orbits along with the local rotation) $\bar{\mathcal{G}'}$ of $\mathcal{O}_F$ such that $\Gamma(\mathcal{O}_F/\langle \bar{\mathcal{G}'}\rangle=(g';n_1',n_2',\dots,n_l')$, where
\begin{center}
$n_i'\in \{n_1,n_2,\dots, n_{\ell}\}\cup \{m_i:\gcd(m_i,n)=1\} \cup \{n_jm_i:\gcd(m_i,n)=1\}\cup \{nm_j\}$
\end{center}
for all $1\leq i\leq l$.
\item $D_G=(mn,g',r'';(c_1',n_1'),(c_1',n_1'),\dots,(c_l',n_l'))$, where
\[
c_i'=
\begin{cases}
c_j, & \text{ if } n_i'=n_j \text{ and }\\
c_j \pmod{n_j}, & \text{ if }n_i'=n_jm_i
\end{cases}
\]
for all $1\leq i\leq l$.
\end{enumerate}
\end{prop} 

\begin{lem}
\label{lem:surface_orbit_power}
Let $S_g(p)$ be the disjoint union of $p$-copies of $S_g$. For $F\in \map(S_g)$ and the $p$-cycle $\sigma_p$, $(\sigma_{p} F)^u$ decomposes into $\gcd(p,u)$-copies of $\sigma_{p/\gcd(p,u)} F^{u/\gcd(p,u)}$.
\end{lem}
\begin{proof}
Write $u=pq+r$, where $0\leq r<p$. We divide the proof into two cases. First consider the case when $r=0$. As in Remark~\ref{rem:cyclic_multiple}, we have
\[
(\sigma_p F)^u=(\sigma_p F^q\sigma_p^{-1})(\sigma_p^2 F^q\sigma_p^{-2})\cdots (\sigma_p^{p-1} F^q\sigma_p^{1-p})F^q.
\]
Since $\gcd(p,u)=p$, $q=u/p$, and $\sigma_p^iF^q\sigma_p^{-i}$ are all supported on connected surfaces, we get $\gcd(p,u)=p$-copies of $\sigma_{p/\gcd(p,u)}F^{u/\gcd(p,u)}$. Now consider the case when $q=0$ and $u=r<p$. To prove the result, we compare the first return map on each connected surface. Denote $\gcd(p,r)$ by $d$. As in Remark~\ref{rem:cyclic_multiple}, we have
\[
(\sigma_{p/d}F^{r/d})^{p/d}=(\sigma_{p/d} F^{r/d}\sigma_{p/d}^{-1})\cdots (\sigma_{p/d}^{(p/d)-1} F^{r/d}\sigma_{p/d}^{1-(p/d)})F^{r/d}.
\]
On the other hand, we have
\[
((\sigma_pF)^r)^{p/d}=((\sigma_pF)^p)^{r/d}=(\sigma_p F^{r/d}\sigma_p^{-1})\cdots (\sigma_p^{p-1} F^{r/d}\sigma_p^{1-p})F^{r/d}
\]
Since $|\sigma_p^r|=p/d$, we have that $(\sigma_pF)^r$ decomposes into $d$-copies of $\sigma_{p/d}F^{r/d}$. In general, since $\gcd(p,u)=\gcd(p,r)$, the result follows.
\end{proof}

\noindent The following result is an immediate consequence of Lemma~\ref{lem:surface_orbit_power}.
\begin{lem}
\label{lem:root_product_surface_orbit}
Let $F=\textstyle\prod_{i'=1}^{s'}(\sigma_{p_{i'}'} F_{i'})$ and $G=\textstyle\prod_{i=1}^s (\sigma_{p_i} G_i)$ be two mapping classes of disconnected surfaces. Then $G^m=F$ if and only if for each $i\in \{1,2,\dots,s\}$, there exists $i'\in \{1,2,\dots,s'\}$ and a positive divisor $u_i$ of $m$ such that
\begin{enumerate}[(i)]
\item $p_{i'}'=p_i/\gcd(u_i,p_i)$.
\item $F_{i'}=G_i^{u_i/\gcd(u_i,p_i)}$.
\item $\sum_{i=1}^{s}p_i=\sum_{i'=1}^{s'}p_{i'}'$ and $\sum_{i=1}^s \gcd(p_i,u_i)=s'$.
\end{enumerate}
\end{lem}

\noindent We will now derive equivalent conditions for the primitivity of pseudo-periodic mapping classes. 

\begin{theorem}
\label{thm:pp_primitive}
For $g \geq 2$, let $F\in \map(S_g)$ be a pseudo-periodic mapping class with
\begin{align*}
\mathcal{D}_F=&\llbracket n';(D_1',p_1'),\dots,(D_{s'}', p_{s'}');\llparenthesis m_1',q_1',k_1';i_1',j_1'\rrparenthesis,\dots,\llparenthesis m_{\ell_1'}',q_{\ell_1'}',k_{\ell_1'}';i_{\ell_1'}', j_{\ell_1'}'\rrparenthesis,\llparenthesis m_{\ell_1'+1}',q_{\ell_1'+1}',k_{\ell_1'+1}';\\&i_{\ell_1'+1}',j_{\ell_1'+1}'\rrparenthesis^*,\dots,\llparenthesis m_{\ell_2'}',q_{\ell_2'}',k_{\ell_2'}';i_{\ell_2'}',j_{\ell_2'}'\rrparenthesis^*;(m_{\ell_2'+1}',q_{\ell_2'+1}',k_{\ell_2'+1}';  i_{\ell_2'+1}',j_{\ell_2'+1}'), \dots,(m_{\ell_3'}',q_{\ell_3'}',k_{\ell_3'}';\\&i_{\ell_3'}', j_{\ell_3'}'),(m_{\ell_3'+1}',q_{\ell_3'+1}',k_{\ell_3'+1}';i_{\ell_3'+1}',j_{\ell_3'+1}')^*,\dots,(m_{\ell'}',q_{\ell'}',k_{\ell'}';i_{\ell'}',j_{\ell'}')^*\rrbracket,
\end{align*}
Then $F$ has a root $G$ of degree $m$ with 
\begin{align*}
\mathcal{D}_G=&\llbracket n;(D_1,p_1),\dots,(D_s, p_s);\llparenthesis m_1,q_1,k_1;i_1,j_1\rrparenthesis,\dots,\llparenthesis m_{\ell_1},q_{\ell_1},k_{\ell_1};i_{\ell_1},j_{\ell_1}\rrparenthesis,\llparenthesis m_{\ell_1+1},q_{\ell_1+1},k_{\ell_1+1};\\&i_{\ell_1+1},j_{\ell_1+1}\rrparenthesis^*,\dots,\llparenthesis m_{\ell_2},q_{\ell_2},k_{\ell_2};i_{\ell_2},j_{\ell_2}\rrparenthesis^*;(m_{\ell_2+1},q_{\ell_2+1},k_{\ell_2+1};i_{\ell_2+1},j_{\ell_2+1}),\dots,(m_{\ell_3},q_{\ell_3},k_{\ell_3};\\&i_{\ell_3},j_{\ell_3}),(m_{\ell_3+1},q_{\ell_3+1},k_{\ell_3+1};i_{\ell_3+1},j_{\ell_3+1})^*,\dots,(m_{\ell},q_{\ell},k_{\ell};i_{\ell},j_{\ell})^*\rrbracket,
\end{align*}
if and only if following conditions hold.
\begin{enumerate}[(i)]
\item For each $i\in \{1,2,\dots,s\}$, there exist $i'\in \{1,2,\dots,s'\}$ and a positive divisor $u_i$ of $m$ such that $p_i=p_{i'}'\gcd(p_i,u_i)$ and $D_i$ carries a root of $D_{i'}'$ of degree $u_i/\gcd(p_i,u_i)$.
\item $n=mn'$ and $\{q_1,q_2,\dots,q_{\ell}\}=\{q_1',q_2'\dots,q_{\ell'}'\}$.
\item $\sum_{i=1}^{s}p_i=\sum_{i'=1}^{s'}p_{i'}'$, $\sum_{r=1}^{\ell}m_r=\sum_{r'=1}^{\ell'}m_{r'}'$, and $\sum_{i=1}^s \gcd(p_i,u_i)=s'$.
\end{enumerate}
\begin{proof}
A cyclic data set $D$ induces a $C_{p\cdot n(D)}$-action $\sigma_{p} D$ on the surface $S_{g(D)}(p)$, where $\sigma_p$ is a cyclic permutation of components of $S_{g(D)}(p)$ (see Remark \ref{rem:cyclic_multiple}). The cyclic action $(\sigma_{p} D)^u$ decomposes into $\gcd(p,u)$-copies of the action $\sigma_{p/\gcd(p,u)} D'$, where $D$ carries a root of $D'$ of degree $u/\gcd(p,u)$. Hence, a root of $F$ can be constructed from roots of its periodic components.  If $\mathcal{D}_G$ carries a root of $\mathcal{D}_F$ of degree $m$, then for each $i\in \{1,2,\dots,s\}$ there exist $i'\in \{1,2,\dots,s'\}$ and a positive divisor $u_i$ of $m$ such that the cyclic action $(\sigma_{p_i}D_i)^{u_i}$ decomposes into $\gcd(p_i,u_i)$-copies of the cyclic action $(\sigma_{p_{i'}'} D_{i'}')$, where $D_i$ carries a root of $D_{i'}'$ of degree $u_i/\gcd(p_i,u_i)$ and $p_{i'}'=p_i/\gcd(p_i,u_i)$. Since $\mathcal{D}_G$ carries a root of $\mathcal{D}_F$ of degree $m$ and correspond to roots of the same multitwist about $\mathcal{C}(F)$, we have $n=mn'$ and $\{q_1,q_2,\dots,q_{\ell}\}=\{q_1',q_2',\dots,q_{\ell'}'\}$. Other necessary conditions follow from a simple counting argument by tallying the number of orbits induced by the actions of $F$ and $G$ on $S_g(\mathcal{C}(F))$ and $\mathcal{C}(F)$. Conversely, by reversing our arguments, we can see that if $\mathcal{D}_F$ and $\mathcal{D}_G$ satisfy the given conditions, then $\mathcal{D}_G$ carries a root of $\mathcal{D}_F$ of degree $m$.
\end{proof}
\end{theorem}

\noindent All possible roots of pseudo-periodic mapping classes can be constructed up to conjugacy by using Theorem \ref{thm:pp_primitive}.

\subsection{Primitivity of arbitrary mapping classes}
For $g\geq 2$, let $F\in \map(S_g)$ be a reducible mapping class. From the Nielsen-Thurston classification \cite{thurston2}, it follows that
\begin{equation}
\label{eqn:canonical_form}
F=G\textstyle\prod_{i=1}^r F_i,
\end{equation}
where $G$ is a pseudo-periodic mapping class, and each $F_i$ is a pseudo-Anosov mapping class on a surface orbit. It follows that $F$ is primitive if and only if at least one of $G$ or $\textstyle\prod_{i=1}^r F_i$ is primitive, and a root of $F$ can be constructed as a product of roots of $G$ and $\textstyle\prod_{i=1}^r F_i$, where the degree of the root will be the least common multiple of individual degrees. Given an arbitrary mapping class $F\in \map(S_g)$ expressed as a product of Dehn twists, we will give an algorithm that solves the general primitivity problem in $\map(S_g)$.

\begin{algo}
\label{algo:gen_prim}
For $g\geq 2$, let $F\in \map(S_g)$ be a mapping class given as a word in Lickorish generators.
\begin{enumerate}[Step 1:]
\item Determining the type of $F$ (i.e. periodic, reducible or pseudo-Anosov) using the curver algorithm.
\item If $F$ is periodic, then proceed to Step \ref{step:periodic}. If $F$ is pseudo-Anosov, then proceed to Step \ref{step:anosov}. If $F$ is infinite order reducible, then proceed to Step \ref{step:reducible}.
\item Construct the cyclic data set (see Remark \ref{rem:mapping_class_to_conj_class}) $D_F$ and use Proposition \ref{prop:periodic_primitive} to determine the primitivity of $F$ and computation of conjugacy classes of roots, if they exists.
\label{step:periodic}
\item Use the flipper algorithm to determine the primitivity of $F$ and computation of conjugacy classes of roots, if they exist. 
\label{step:anosov}
\item Compute $\mathcal{C}(F)$ using the curver algorithm. Write the mapping class $F$ in the canonical form as in Equation (\ref{eqn:canonical_form}).
\begin{enumerate}[Step 5a:]
\item Determine the primitivity of pseudo-periodic component from Step \ref{step:pseudo-periodic}.
\item Determine the primitivity of $\prod_{i=1}^r F_i$ from Lemma~\ref{lem:root_product_surface_orbit} and Step \ref{step:anosov}.
\end{enumerate}
\label{step:reducible}
\item Construct the pseudo-periodic data set $\mathcal{D}_G$ (see Remark~\ref{rem:mapping_class_to_conj_class}) and use Theorem \ref{thm:pp_primitive} to determine the primitivity of $\mathcal{D}_G$ and computation of conjugacy classes of roots, if they exist.
\label{step:pseudo-periodic}
\end{enumerate}
\end{algo}

\begin{rem}
\label{rem:mapping_class_to_conj_class}
Given a pseudo-periodic mapping class $F$ written as a word in Dehn twists, it is possible to determine $\mathcal{D}_F$. This process will require an appropriate application of Theorem~\ref{thm:main_theorem}, the curver algorithm \cite{curver}, and an understanding of the orbits of carefully chosen curves (that either lie in the reduction systems or outside).

We will now try to explain the process of extracting a data set $D_F$ from a given pseudo-periodic $F$. Let $F\in \map(S_g)$ be a pseudo-periodic mapping class given as a product of Lickorish Dehn twist. First we determine the periodic components of $F$ and $\mathcal{C}(F)$ using the curver algorithm. The conjugacy invariants of periodic components of $F$, that is, their cyclic data set can be obtained by the curver algorithm and analyzing the size of orbits of suitably chosen curves lying either in $\mathcal{C}(F)$ or outside. By taking $F^n$, where $n$ is the least common multiple of orders of periodic components of $F$, the exponents of Dehn twists about curves in the $\mathcal{C}(F)$ can be obtained. The twist factor associated to the orbit of curves of $\mathcal{C}(F)$ can be determined by taking sum of induced angles about compatible cone points as in Definition \ref{defn:twist_factor}. Consequently, the pseudo-periodic data set $D_F$ corresponding to $F$ can be determined.
\end{rem}

\noindent If a mapping class is non-primitive, then an upper bound on the degree of roots can be derived using Proposition \ref{prop:general_bounds} and Lemma \ref{lem:pa_bound}.

\begin{exmp}
For simple closed curves shown in Figure \ref{fig6}, let $F=(T_{a_1}T_{b_1})^3(T_{a_2}T_{b_2})^4T_c^{-1}\in \map(S_2)$. Since $(T_{a_1}T_{b_1})^6=T_c=(T_{a_2}T_{b_2})^6$, we have $F^6=T_c$. Hence, $F$ is a pseudo-periodic mapping class with $\mathcal{C}(F)=\{c\}$. The periodic components of $F$ are $(T_{a_1}T_{b_1})^3$ and $(T_{a_2}T_{b_2})^4$ in $\map(S_1)$.
\begin{figure}[ht]
\labellist
\small
\pinlabel $b_2$ at 180, 85
\pinlabel $b_1$ at 65, 85
\pinlabel $a_2$ at 177, -8
\pinlabel $a_1$ at 68, -8
\pinlabel $c$ at 118, 15
\endlabellist
\centering
\includegraphics[scale=.8]{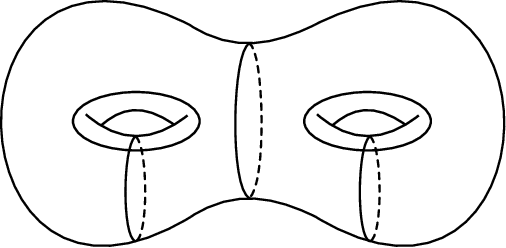}
\caption{Simple closed curves on $S_2$.}
\label{fig6}
\end{figure}
\noindent The pseudo-periodic data set corresponding to conjugacy class of $F$ is given as
\begin{center}
$\mathcal{D'}=\llbracket 6;(D_1,1),(D_2,1);\llparenthesis 1,1,1;1,2\rrparenthesis;-\rrbracket$,
\end{center}
where $D_1=(2,0;(1,2),(1,2),(1,2),(1,2)_1)$, $D_2=(3,0;(2,3),(2,3),(2,3)_1)$, and the cone points with the same suffix are compatible with twist factor $1$. The mapping class $F$ is not primitive and has a root of degree $2$ given by the pseudo-periodic data set
\begin{center}
$\mathcal{D}=\llbracket 12;(E_1,1),(E_2,1);\llparenthesis 1,1,1;1,2\rrparenthesis;-\rrbracket$,
\end{center}
where $E_1=(4,0;(1,2),(1,4),(1,4)_1)$, $E_2=(6,0;(1,2),(2,3),(5,6)_1)$ and the cone points with the same suffix are compatible with twist factor $1$. For $i=1,2$, we observe that $E_i$ carries a root of $D_i$ of degree $2$. As an application of Proposition \ref{prop:periodic_primitive} and Theorem \ref{thm:pp_primitive} we have that $D$ is the only conjugacy class which carries a root of $D'$. For the mapping class $G=T_{a_1}^2T_{b_1}(T_{a_2}T_{b_2})^5T_c^{-1}\in \map(S_2)$, we have $\mathcal{D}_G$ and $\mathcal{D}$ are equivalent. Further, a direct computation shows that $G^2=F$.   
\end{exmp}

\section{Realizable bounds on the primitivity of pseudo-periodic mapping classes}
\label{sec:bounds}
In this section, we apply Proposition~\ref{thm:main_theorem} to derive realizable bounds on the degrees of roots of pseudo-periodics in $\map(S_g)$, the Torelli group $\mathcal{I}(S_g)$, the level-$m$ subgroup $\map(S_g)[m]$, and the commutator subgroup of $\map(S_2)$. For deriving these bounds, we will make extensive use of the Riemann-Hurwitz equation that appeared in condition (vi) of Definition \ref{defn:cyclic_data}. For a $C_n$-action on a surface $S_g$ with quotient orbifold $\mathcal{O}=S_g/C_n$ of signature $(g_0;n_1,n_2,\dots,n_k)$, the Riemann-Hurwitz equation is given by:
\begin{center}
$\frac{2g-2}{n}=2g_0-2+\sum_{i=1}^{k}\left (1-\frac{1}{n_i}\right)$.
\end{center}
\noindent For the realization of upper bounds in the following and subsequent propositions, we need following technical lemma from \cite[Lemma 8.1]{rajeevsarathy3}.

\begin{lem}
\label{lem:gcd_lemma}
Let $d_1$, $d_2$ be relatively prime integers, an let $Q$ be a finite set of primes. If $2\in Q$, assume that $d_1$ and $d_2$ are not both odd. Then there exist integers $c_1$ and $c_2$ such that $c_1d_1+c_2d_2=1$ and neither $c_1$ nor $c_2$ is divisible by any primes in $Q$.
\end{lem}

\subsection{Realizable bounds on the degree in the mapping class group}
Recall that there exists a unique positive integer $q(F)=q(\mathcal{D}_F)$ associated to the pseudo-periodic mapping class $F$.

\begin{prop}
\label{prop:general_bounds}
For $g\geq 2$, let $F\in \map(S_g)$ be a pseudo-periodic mapping class with canonical reduction system $\mathcal{C}(F)$. Assume that $F$ has a root $G$ of degree $n$.
\begin{enumerate}[(i)]
\item If $\mathcal{C}(F)$ is nonseparating, then $3 \leq n\leq q(F)(2g-1)$. These bounds are always realized by roots of powers of a Dehn twist about nonseparating curve.
\item If $\mathcal{C}(F)$ is separating, then $2\leq n\leq 3q(F)(g+1)(g+2)$. This upper bound is only realized by a root of powers of a Dehn twist about separating curve of genus $[g/2]$, where $g\equiv 0,9 \pmod{12}$.
\end{enumerate}
\begin{proof}
Since a nontrivial pseudo-periodic is a root of a multitwist, it suffices to assume that $F$ is a multitwist. We first consider the case when canonical reduction system $\mathcal{C}(F)$ is nonseparating. We assume that $|\mathcal{C}(F)|=k$. The pseudo-periodic  data set corresponding to $G$ can be written as
\begin{align*}
\mathcal{D}_G=\llbracket n;(D_1,1);-;(m_1,q_1,k_1;1,1),\dots,
&(m_{\ell_1},q_{\ell_1},k_{\ell_1};1,1),\\
&(m_{\ell_1+1},q_{\ell_1+1},k_{\ell_1+1};1,1)^*,\dots,(m_{\ell},q_{\ell},k_{\ell};1,1)^*\rrbracket.
\end{align*}
\noindent We have $n=(q_r/k_r)n(D_1)$ for every $1\leq r\ell$. We write $n_1$ for $n(D_1)$. If there exist some $r'$ such that $m_{r'}=n_1$, then $n_1=m_{r'}\leq k\leq g$. Therefore, we assume that $m_r<n_1$ for every $r$. First consider the case when there is a single orbit of size $1$, that is, $k=1$. If $G$ preserves the sides of the curve, then it has been shown in \cite[Corollary 2.6, 2.9]{rajeevsarathy4} that 
\[
3\leq n_1 \leq 
\begin{cases}
2g-1, & \text{if } \gcd(k_1,n_1)=1,\\
4(g-1), & \text{ otherwise.}
\end{cases}
\]
\noindent When $\gcd(k_1,n_1)>1$, we must have $k_1\geq 2$. In the case when $G$ exchanges the sides of the curve, it was shown in \cite[Remark 3.6]{rajeevsarathy4} that $n_1\leq 4g-2$. It was shown in \cite{rajeevsarathy3} that Dehn twist can not have side-exchanging roots, therefore $k_1\geq 2$.  Now we assume that $k>1$ and $G$ preserves sides of all the curves. An orbit of size $m_r$ corresponds to two cone points of order $n/m_r$ in the cyclic data set $D_1$. If $\ell>1$, then from Riemann-Hurwitz equation, it follows that
$$n_1\leq \frac{g-1}{\ell-1}\leq g-1.$$
When $\ell=1$ we observe that $G^k$ fixes all the curves of $\mathcal{C}(F)$ and therefore, from Riemann-Hurwitz equation, we have
$$\frac{n_1}{k} \leq \frac{g-1}{k-1},$$
and hence, $n_1\leq 2(g-1)$. From Equation (\ref{eqn:twist_compa}), it follows that if $m_r>1$, then $m_r\mid k_r$, therefore, $k_r\geq 2$. Considering all the cases described above, we have
\begin{equation}
\label{eqn:nonsep_bound}
3\leq n\leq q(F)(2g-1).
\end{equation}

Now we consider the case when the canonical reduction system $\mathcal{C}(F)$ is separating. We first consider the case when $\mathcal{C}(F)$ is a bounding multicurve. For $i=1,2$, we write $n_i$ and $g_i$ for $n(D_i)$ and $g(D_i)$, respectively for simplicity. If $G$ preserves the sides of curves of $\mathcal{C}(F)$, then the pseudo-periodic data set corresponding to $G$ can be written as
\begin{align*}
\mathcal{D}_G=\llbracket n;(D_1,1),(D_2, 1);\llparenthesis m_1,q_1,k_1;1,2\rrparenthesis,\dots,\llparenthesis m_{\ell},q_{\ell},k_{\ell};1,2\rrparenthesis;-\rrbracket.
\end{align*}
For every $r$, we have that $n=(q_r/k_r)\lcm(n_1,n_2)$, and therefore, $n\leq q(F)n_1n_2$. First assume that $k=1$, that is, $\mathcal{C}(F)$ contains a single separating curve $c$. Each $D_i$ has a distinguished fixed point obtained by $c$ capping boundary components of $S_g\setminus c$. It is well-known \cite{harvey} that $n_i\leq 4g_i+2$. When $n_i$ is odd, it follows from the Riemann-Hurwitz equation that, $n_i\leq 3g_i+3$. If both of $n_i$'s are even, then $$n\leq 2q(F)(4g_1g_2+2g+1).$$ If at least one of $n_i$ is odd (say $n_1$) then $$n\leq 6q(F)(2g_1g_2+2g_1+g_2+1).$$ It can be easily shown that $(2g_1g_2+2g_1+g_2+1)$ attains its maximum value at $g_1=g_2=g/2$ when $g$ is even and at $g_1=(g+1)/2$, $g_2=(g-1)/2$ when $g$ is odd subject to the conditions that $g_1+g_2=g$ and $g_i$ is a positive integer. Hence,
$$2 \leq n\leq 3q(F)(g+1)(g+2).$$
\noindent Using the Riemann-Hurwitz equation and arguing as before, it can be established that
\[
n\leq
\begin{cases}
q(F)(g^2-1), & \text{if } \ell=1, k>1 \text{ and}\\
3q(F)g(g-2), & \text{if } \ell>1.
\end{cases}
\]
\noindent If $G$ exchanges the sides of curves of $\mathcal{C}(F)$, then the pseudo-periodic data set corresponding to $G$ can be written as 
\begin{align*}
\mathcal{D}_G=\llbracket n;(D_1,2);\llparenthesis 1,q_1,k_1;1,1\rrparenthesis^*,\dots,\llparenthesis 1,q_m,k_m;1,1\rrparenthesis^*;-\rrbracket.
\end{align*}
Since $n\leq 2q(F)n_1$ and $g_1=g/2$, we have that
$$n\leq 4q(F)(2g_1+1)=4q(F)(g+1).$$
In general, assume that the canonical reduction system $\mathcal{C}(F)$ is not a bounding multicurve. Let $\mathcal{D}_G$ be the pseudo-periodic data set corresponding to $G$ as in Definition \ref{defn:encoding}. For each $r$, we have
$$n=\left( \frac{q_r}{k_r}\right) p_{i_r}p_{j_r}\lcm \left(\frac{n(D_{i_r})}{p_{j_r}},\frac{n(D_{j_r})}{p_{i_r}} \right).$$
As shown before, we have
\begin{equation}
\label{eqn:sep_bound}
n\leq 3q(F)(g+1)(g+2).
\end{equation}

Now we construct some pseudo-periodic mapping class realizing bounds obtained above. For a positive integer $q$, the upper bound obtained in Equation (\ref{eqn:nonsep_bound}) is realized by the pseudo-periodic with data 
\begin{center}
$\mathcal{D}_G=\llbracket q(2g-1);(D_1,1);-;(1,q,1;1,1)\rrbracket$,
\end{center}
which corresponds to a conjugacy class of a root of $F=T_c^q$ of degree $q(2g-1)$, where
\begin{center}
$D_1=(2g-1,0;(2g-5,2g-1),(2,2g-1)_1,(2,2g-1)_1)$,
\end{center}
where cone points with the same suffix are compatible with twist factor $1$, and $c$ is a nonseparating curve in $S_g$. For realization of the upper bound obtained in Equation (\ref{eqn:sep_bound}), we consider the cases when $g$ is even and odd separately.

First assume that $g$ is even. For $g=12s$, take $g_1=g_2=6s$, where $s$ is a positive integer. Consider the cyclic data sets
\begin{align*}
&D_1=(4g_1+2,0;(a_1,2),(a_1g_1,2g_1+1),(a_1,4g_1+2)) \text{ and}\\
&D_2=(3g_2+3,0;(2a_2,3),(a_2(g_2/3),g_2+1),(a_2,3g_2+3)).
\end{align*}
Taking $n_1=4g_1+2$ and $n_2=3g_2+3$, we observe that $\gcd(n_1,n_2)=1$, consequently $\lcm(n_1,n_2)=n_1n_2=3(g+1)(g+2)=n$ (say). Since both $n_i$'s are not odd and $\gcd(n/n_1,n/n_2)\linebreak=1$, from Lemma \ref{lem:gcd_lemma}, it is possible to choose integers $p_1$ and $p_2$ such that $\gcd(p_i,n_i)=1$ and
\begin{center}
$\frac{n}{n_1}p_1+\frac{n}{n_2}p_2\equiv 1\pmod n$.
\end{center}
For $a_i=p_i^{-1}\pmod {n_i}$, we observe that the fixed points of $D_1$ and $D_2$ are compatible with twist factor $1$, and hence, the pseudo-periodic with data
\begin{center}
$\mathcal{D}_G=\llbracket 3q(g+1)(g+2);(D_1,1),(D_2,1);\llparenthesis 1,q,1;1,2\rrparenthesis;-\rrbracket$
\end{center}
corresponds to a conjugacy class of a root of $F=T_c^q$ of degree $3q(g+1)(g+2)$, where $c$ is a separating curve in $S_g$ of genus $g/2$ and $12\mid g$.

Now, we assume that $g$ is odd. For $g=12s+9$, take $g_1=6s+5$ and $g_2=6s+4$, where $s$ is a positive integer. Consider the cyclic data sets
\begin{align*}
&D_1=(4g_1+2,0;(a_1,2),(a_1g_1,2g_1+1),(a_1,4g_1+2)) \text{ and}\\
&D_2=(3g_2+3,0;(a_2,3),(a_2(2g_2+1)/3,g_2+1),(a_2,3g_2+3)).
\end{align*}
For $n_1=4g_1+2$ and $n_2=3g_2+3$, since $\gcd(n_1,n_2)=1$, consequently $\lcm(n_1,n_2)=n_1n_2=3(g+1)(g+2)=n$ (say). Since both $n_i$'s are not odd and $\gcd(n/n_1,n/n_2)=1$, from Lemma \ref{lem:gcd_lemma}, it is possible to choose integers $p_1$ and $p_2$ such that $\gcd(p_i,n_i)=1$ and
\begin{center}
$\frac{n}{n_1}p_1+\frac{n}{n_2}p_2\equiv 1\pmod n$.
\end{center}
For $a_i=p_i^{-1}\pmod {n_i}$, we observe that the fixed points of $D_1$ and $D_2$ are compatible with twist factor $1$, and hence, the pseudo-periodic with data
\begin{center}
$\mathcal{D}_G=\llbracket 3q(g+1)(g+2);(D_1,1),(D_2,1);\llparenthesis 1,q,1;1,2\rrparenthesis;-\rrbracket$
\end{center}
corresponds to a conjugacy class of a root of $F=T_c^q$ of degree $3q(g+1)(g+2)$, where $c$ is a separating curve in $S_g$ of genus $(g-1)/2$ and $g\equiv 9 \pmod{12}$.
\end{proof}
\end{prop}

\noindent Following corollaries are immediate from Proposition \ref{prop:general_bounds}.

\begin{cor}
The highest degree of a root of Dehn twist about a separating curve in $\map(S_g)$ is $3(g+1)(g+2)$. This upper bound is realized for infinitely many $g$.
\end{cor}

\begin{cor}
\label{cor:bound_components}
The upper bound realizing roots of pseudo-periodic mapping classes in $\map(S_g)$ are primitive pseudo-periodic mapping classes that decompose canonically into irreducible periodic mapping classes whose Nielsen representatives have at least one fixed point. Any primitive pseudo-periodic mapping class must be the highest degree root of a multitwist.
\end{cor}

\subsection{Realizable bounds on the degree in the Torelli group and the level-$m$-subgroups}
This subsection derives realizable bounds on the degrees of roots of pseudo-periodic mapping classes in the Torelli group and the level-$m$-subgroup of $\map(S_g)$. In Proposition \ref{prop:general_bounds}, we have seen that a root of a Dehn twist about a separating curve realizes the highest degree in $\map(S_g)$. Since the powers of a Dehn twist about a separating curve lies in $\map(S_g)[m]$ for every $m\geq 2$, in the following proposition we will assume that there are no separating curves in the canonical reduction system.

\begin{prop}
\label{prop:torelli_levelm}
For $g \geq 3$ and $m \geq 2$, let $F\in \map(S_g)[m]$ be a pseudo-periodic mapping class. Suppose that the canonical reduction system for $F$ does not contain any separating curves and that $F$ has a root in $\map(S_g)$ of degree $n$.
\begin{enumerate}[(i)]
\item If $F\in \mathcal{I}(S_g)$, then $2 \leq n \leq q(F)g(g-2)$. The upper bound for $n$ is realized only when $F = (T_{c_1}T_{c_2}^{-1})^q$, where $\{c_1,c_2\}$ is a bounding pair of genus $(g/2)-1$ with $2\mid g$ and $q$ is a positive integer.
\item If $F\in \map(S_g)[m]\setminus \mathcal{I}(S_g)$, then $m \leq n \leq 3q(F)g(g-2)$. The upper bound for $n$ is realized when $F = (T_{c_1}T_{c_2}^k)^q$, where $\{c_1,c_2\}$ a bounding pair of genus $(g/2)-1$ with $g \equiv -4 \pmod{24}$, $q$ is a positive integer, and $k=(1/4)g(g-2) - 1$.
\end{enumerate}
\begin{proof}
Since a nontrivial pseudo-periodic is a root of a multitwist, it suffices to assume that $F$ is a multitwist. We observe that the canonical reduction system $\mathcal{C}(F)$ is a disjoint union of bounding multicurves and a nonseparating multicurve. In the proof of Proposition \ref{prop:general_bounds}, we have seen that the highest degree of a root is achieved when $\mathcal{C}(F)$ is a bounding pair (as $\mathcal{C}(F)$ does not contain separating curves). Therefore, for deriving bounds on the degrees of roots of $F$, it suffices to consider $F=T_{c_1}^{k_1}T_{c_2}^{k_2}$ and
\begin{center}
$\mathcal{D}_G=\llbracket n;(D_1,1),(D_2,1);\llparenthesis 1,k_1,k_1;1,2\rrparenthesis,\llparenthesis 1,k_2,k_2;1,2\rrparenthesis;-\rrbracket.$
\end{center}
We observe that $F\in \map(S_g)[m]$ (resp. $F\in \mathcal{I}(S_g)$) if and only if $m\mid k_1+k_2$ (resp. $k_1+k_2=0$). For $i=1,2$, we write $n_i$ and $g_i$ for $n(D_i)$ and $g(D_i)$, respectively for simplicity. Since each $D_i$ has at least two distinguished fixed points, from the Riemann-Hurwitz equation, we have that $n_i\leq 4g_i$ when $n_i$ is even, and $n_i\leq 3g_i$ when $n_i$ is odd. Since $n=\lcm(n_1,n_2)$, we have
\begin{equation}
\label{eqn:levelm_bound1}
n=(n_1n_2)/\gcd(n_1,n_2)\leq 12g_1g_2\leq 3g(g-2),
\end{equation}
where $g_1g_2$ attains its maximum value at $g_1=g/2$ and $g_2=g/2-1$ subject to the conditions that $g=g_1+g_2+1$, $g_1\neq g_2$, and $g_i$ is a positive integer.

The upper bound in Equation (\ref{eqn:levelm_bound1}) is realized if and only if $g$ is even, $n_1=2g$ and $n_2=(3/2)(g-2)$, where $\gcd(n_1,n_2)=1$. The periodic components of the bound-realizing roots are given by the data sets
\begin{align*}
&D_1=(n_1,0;(a_1,n_1)_1,(a_2,n_1)_2,(1,2)) \text{ and}\\
&D_2=(n_2,0;(b_1,n_2)_1,(b_2,n_2)_2,(1,3)),
\end{align*}
where $a_1+a_2\equiv g \pmod{n_1}$ and $b_1+b_2\equiv (g-2)\pmod{n_2}$. Compatibility of the cone points with same suffix implies that
\begin{align}
\label{eqn:compatibility_levelk}
(3/2)g(g-2)a_1^{-1}a_2^{-1}+2g(g-2)b_1^{-1}b_2^{-1}\equiv k_1+k_2 \pmod {3g(g-2)}.
\end{align}
From Equation (\ref{eqn:compatibility_levelk}), it follows that if $m\mid (g/2)(g-2)$, then $F\in \map(S_g)[m]$ which has a root of degree $3g(g-2)$. Again, from Equation (\ref{eqn:compatibility_levelk}), it can be seen that if $k_1+k_2=0$, then either $2$ or $3/2$ divide $1$, which is impossible. Arguing as before, it follows that if $\gcd(n_1,n_2)=1$ and either $n_1>g$ or $n_2>(g-2)$, then the bounding pair maps cannot have a root of degree $\lcm(n_1,n_2)$. Now, we find a realizable bound on the degree of a root of a bounding pair map. As seen before, if $\gcd(n_1,n_2)=1$, then $n\leq n_1n_2\leq g(g-2).$ When $\gcd(n_1,n_2)>1$, since $n_i\leq 4g_i$, we have that $n\leq 4g_1g_2\leq g(g-2),$ where $g_1g_2$ attains its maximum value when $(g_1,g_2)=(g/2,g/2-1)$.

In general, let $F\in \map(S_g)[m]$ having a root $G$ of degree $n$ with corresponding pseudo-periodic data set $\mathcal{D}_G$ as in Definition \ref{defn:encoding}. We have
$$n=\left( \frac{q_r}{k_r}\right) p_{i_r}p_{j_r}\lcm \left(\frac{n(D_{i_r})}{p_{j_r}},\frac{n(D_{j_r})}{p_{i_r}} \right).$$
As shown above, if $F\in \mathcal{I}(S_g)$, then
\begin{equation}
\label{eqn:torelli_bound}
2 \leq n \leq q(F)g(g-2),
\end{equation}
and if $F\in \map(S_g)[m]\setminus\mathcal{I}(S_g)$, then
\begin{equation}
\label{eqn:levelm_bound}
m \leq n \leq 3q(F)g(g-2).
\end{equation}

We will now construct roots that realize the upper bounds obtained in Equation (\ref{eqn:torelli_bound}) and (\ref{eqn:levelm_bound}). For $g=2s$, let $g_1=s$ and $g_2=s-1$, where $s$ is a positive integer. Consider the following data sets
\begin{align*}
&D_1=(4s,0;(1,2),(4s-1,4s)_1,(2s+1,4s)_2) \text{ and}\\
&D_2=(4s-4,0;(1,2),(1,4s-4)_1,(2s-3,4s-4)_2).
\end{align*}
We have that $\lcm(4s,4s-4)=g(g-2)$ and the cone points with the same suffix are compatible with twist factor $\pm 1$. Hence, the pseudo-periodic data set
\begin{center}
$\mathcal{D}_G=\llbracket qg(g-2);(D_1,1),(D_2,1);\llparenthesis 1,q,1;1,2\rrparenthesis,\llparenthesis 1,-q,-1;1,2\rrparenthesis;-\rrbracket$
\end{center}
corresponds to a conjugacy class of a root of $F=(T_{c_1}T_{c_2}^{-1})^q$ of degree $qg(g-2)$, where $q=q(F)$. Now, we consider the following cyclic data sets
\begin{align*}
&D_1=(4g_1,0;(a_1,2),(a_1,4g_1)_1,(a_1(2g_1-1),4g_1)_2) \text{ and}\\
&D_2=(3g_2,0;(a_2,3),(a_2,3g_2)_1,(a_2(2g_2-1),3g_2)_2).
\end{align*}
For $n_1=4g_1$ and $n_2=3g_2$, we have $\gcd(n_1,n_2)=1$, and therefore
\begin{center}
$\lcm(n_1,n_2)=12s(s-1)=3g(g-2).$
\end{center}
For an integer $r\geq 0$ taking $s=12r+10$, we choose integers $a_i$ such that $a_1^{-1}=4r+3=(s-1)/3$ and $a_2^{-1}=-(3r+2)=-(s-2)/4$. We observe that $(2g_1-1)^{-1}\equiv 2g_1-1 \pmod{n_1}$ and $(2g_2-1)^{-1}=g_2-1 \pmod{n_2}$. The cone points with the same suffix are compatible with twist factor $1$ and $(1/4)g(g-2)-1=k$ (say). Hence, the pseudo-periodic data set
\begin{center}
$\mathcal{D}_G=\llbracket 3qg(g-2);(D_1,1),(D_2,1);\llparenthesis 1,q,1;1,2\rrparenthesis,\llparenthesis 1,kq,k;1,2\rrparenthesis;-\rrbracket$
\end{center}
corresponds to a root of $F=(T_{c_1}T_{c_2}^k)^q$ of degree $3qg(g-2)$, where $q=q(F)$.
\end{proof}
\end{prop}

\begin{cor}
For a bounding pair $\{c_1,c_2\}$ in $S_g$, the highest degree of a root of a bounding pair map $T_{c_1}T_{c_2}^{-1}$ is $g(g-2)$. This upper bound is realized for infinitely many $g$.
\end{cor}

\noindent For $m\geq 3$, it is known \cite[Corollary 1.8]{ivanov} that any pseudo-periodic mapping class in $\map(S_g)[m]$ is a multitwist. The following corollary now follows from the Equation (\ref{eqn:compatibility_levelk}) described in the proof of the Proposition \ref{prop:torelli_levelm}. We note that a power of a Dehn twist about a separating curve is always non-primitive.

\begin{cor}
For $g,m\geq 3$, let $F\in \map(S_g)[m]\setminus \mathcal{I}(S_g)$ be a pseudo-periodic mapping class. Suppose that the canonical reduction system for $F$ does not contain any separating curves and $F$ is not a power of a multitwist in $\map(S_g)$. If $m>3g(g-2)$, then $F$ is primitive.
\end{cor}

\subsection{Realizable bounds on the degrees of roots in the commutator subgroup of $\map(S_2)$}
For $g\geq 3$, it is known \cite{harer, powell} that $\map(S_g)$ is perfect. A surjective homomorphism from the braid group $B_6$ onto $\map(S_2)$ can be constructed, which maps the derived series of $B_6$ onto the derived series of $\map(S_2)$. For $n\geq 5$, it is known \cite{lin} that the commutator subgroup of $B_n$ is perfect. We denote the commutator subgroup of $\map(S_2)$ by $\map(S_2)^{(1)}$. It follows that $\map(S_2)\triangleright \map(S_2)^{(1)}$ is the derived series of $\map(S_2)$.

\begin{prop}
\label{prop:bound_derived}
Let $F=T_{c_1}^{q_1}T_{c_2}^{q_2}T_{c_3}^{q_3}\in \map(S_2)$ be a multitwist, where $c_1$ and $c_3$ are nonseparating curves. Then $F\in \map(S_2)^{(1)}$ if and only if either $q_1+q_2+q_3\equiv 0 \pmod {10}$ (resp. $q_1+2q_2+q_3\equiv 0 \pmod{10}$) depending upon whether $c_2$ is a nonseparating (resp. separating) curve. Now we assume that $F\in \map(S_2)^{(1)}$ has a root in $\map(S_2)$ of degree $n$. Then
\[
2q(F)\leq n\leq 
\begin{cases}
3q(F), & c_2 \text{ is nonseparating},\\
12q(F), & c_2 \text{ is separating}.
\end{cases}
\]
Further, these bounds are realized when $F=T_c^q$, where $q$ is a positive integer such that $5\mid q$ and $c$ is a separating curve. The bounds realizing roots decompose canonically into irreducible periodic mapping classes whose Nielsen representatives have at least one fixed point.
\begin{proof}
It is known \cite{birman2} that $\map(S_2)/\map(S_2)^{(1)}\cong \mathbb{Z}_{10}$. Let $\pi:\map(S_2)\to \mathbb{Z}_{10}$ be the natural projection onto the abelianization of $\map(S_2)$. The map $\pi$ maps a nonseparating Dehn twist to $1$ and a separating Dehn twist to $2$. Now it follows that $F\in \map(S_2)^{(1)}$ if and only if either $q_1+q_2+q_3\equiv 0 \pmod {10}$ (resp. $q_1+2q_2+q_3\equiv 0 \pmod{10}$) depending upon whether $c_2$ is a nonseparating (resp. separating) curve. The verification of the realizable bounds is straightforward as the least common multiple of orders of periodic components induced by the root on the surface $S_2(\mathcal{C}(F))$ can be at most $12$.  
\end{proof}
\end{prop}

\subsection{An upper bound on the degree of a root of a pseudo-Anosov} The following lemma which is included simply for the sake of completion, is a rather straightforward application of a result of Penner \cite{penner}.

\begin{lem}
\label{lem:pa_bound}
For $g\geq 2$, let $F\in \map(S_g)$ be a pseudo-Anosov mapping class with stretch factor $\lambda$ which has a root of degree $n$. Then
\begin{center}
$n \leq \frac{12(g-1)\log{\lambda}}{\log2}.$
\end{center}
\begin{proof}
The stretch factor of the root of $F$ must be $\lambda^{\frac{1}{n}}$, and therefore, $\log{\lambda}\geq n \log{\lambda_g}$, where $\lambda_g$ denotes the smallest stretch factor for a pseudo-Anosov in $\map(S_g)$. By \cite{penner}
\begin{center}
$\log{\lambda_g}\geq \frac{\log2}{12(g-1)}$
\end{center}
from which our assertion follows.
\end{proof}
\end{lem}

\section{Normal closure of pseudo-periodic mapping classes}
\label{sec:normal_gen}
This section discusses the normal generation of $\map(S_g)$ by a single pseudo-periodic mapping class. Recently, Margalit and Lanier \cite[Theorem 1.1]{normal} have proved that when $g\geq3$, any periodic mapping class, which is not a hyperelliptic involution, normally generates $\map(S_g)$. For this purpose, we will use the main results from \cite{rajeevsarathy1}, where the following categorization of cyclic actions was introduced.
\begin{defn}
\label{defn:type}
Let $D$ be a cyclic action on $S_g$. Then we say that $D$ is
\begin{enumerate}[(i)]
\item \textit{a rotational action} if either $r(D)\neq 0$ or $D$ is of the form
$$(n,g_0;\underbrace{(c,n),(n-c,n),\dots,(c,n),(n-c,n)}_{\ell-\text{pairs}}),$$
where $\ell \geq 1$, $0<c<n$, $\gcd(c,n)=1$, and $\ell=1$ if and only if $n>2$.
\item \textit{of Type 1}, if $k=3$ and $n_j=n$ for some $j$.
\item \textit{of Type 2}, if $D$ is neither rotational nor of Type 1.
\end{enumerate}
\end{defn}

\begin{lem}[{\cite[Lemmas 2.2-2.4]{normal}}]
\label{lem:wellsuited_curve}
For $g\geq 3$, let $F\in \map(S_g)$, $c$ be a nonseparating curve, and $d$ be a separating curve in $S_g$. The geometric intersection number between isotopy classes of two simple closed curves $c_1$ and $c_2$ will be denoted by $i(c_1,c_2)$. If one of the following conditions hold, then $F$ normally generates $\map(S_g)$. 
\begin{enumerate}[(i)]
\item $i(c,F(c))=1$.
\item $i(c,F(c))=0$ and $[c]\neq [F(c)]$ in homology (ignoring the orientation of curves).
\item $i(d,F(d))\leq 2$.
\end{enumerate}
\end{lem}

\begin{prop}
\label{prop:pp_normal}
For $g\geq 3$, let $F\in \map(S_g)$ be a pseudo-periodic mapping class with at least one nontrivial periodic component $F'$ that is not a hyperelliptic involution. Then $\map(S_g)$ is normally generated by $F$.
\begin{proof}
Let $S$ be the $\langle F\rangle$-invariant subsurface $S_{g(D_{F'})}$ of $S_g$ supporting $F'$. If $D_{F'}$ is a rotational cyclic action, then it is always possible to find a nonseparating curve $c$ in $S$ such that $i(c,F'(c))=0$ and $[c]\neq [F(c)]$ except in the case of a hyperelliptic involution. Figure \ref{fig:rotational} illustrates the existence of such a curve $c$ for rotational actions on $S_3$ and $S_4$. Since $S$ is an $\langle F\rangle$-invariant subsurface, $i(c,F'(c))=i(c,F(c))$. The assertion now follows from Lemma \ref{lem:wellsuited_curve}. 
\begin{figure}[H]
\begin{subfigure}{0.32\textwidth}
\labellist
\small
\pinlabel $\frac{2\pi}{3}$ at 60, 255
\endlabellist
\centering
\includegraphics[scale=.29]{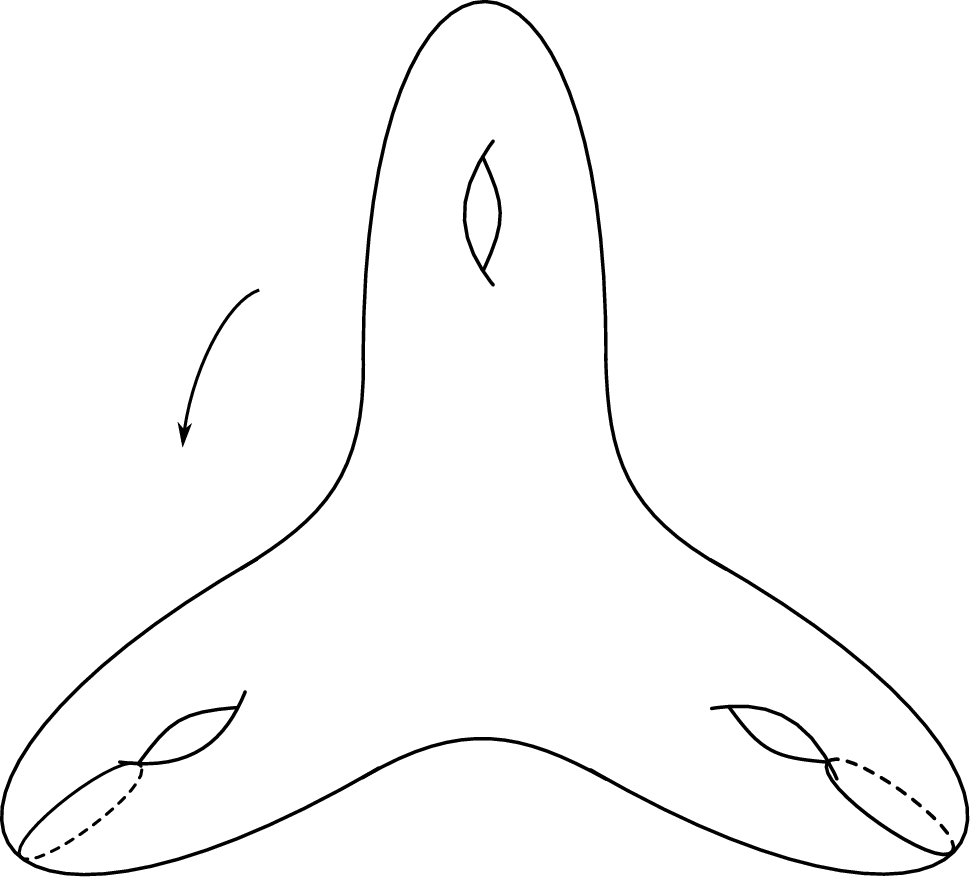}
\end{subfigure}
\hfill
\begin{subfigure}{0.32\textwidth}
\labellist
\small
\pinlabel $\frac{2\pi}{3}$ at 60, 255
\endlabellist
\centering
\includegraphics[scale=.29]{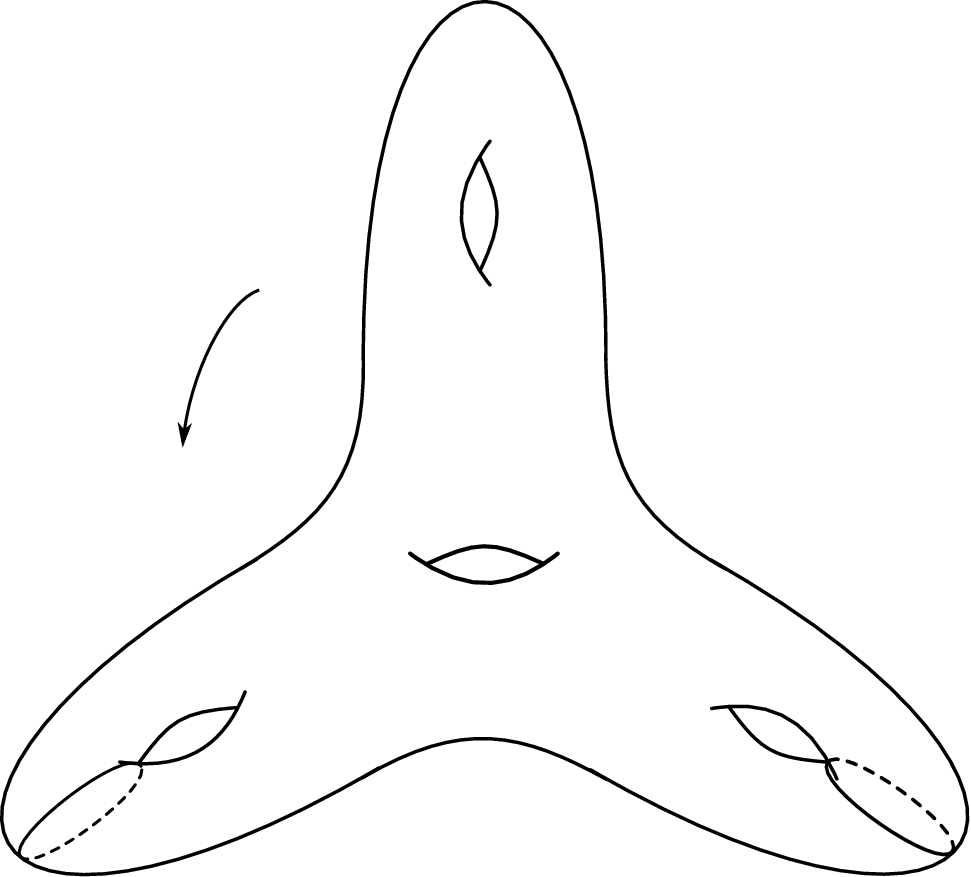}
\end{subfigure}
\hfill
\begin{subfigure}{0.32\textwidth}
\labellist
\small
\pinlabel $\pi$ at 415, 225
\endlabellist
\centering
\includegraphics[scale=.29]{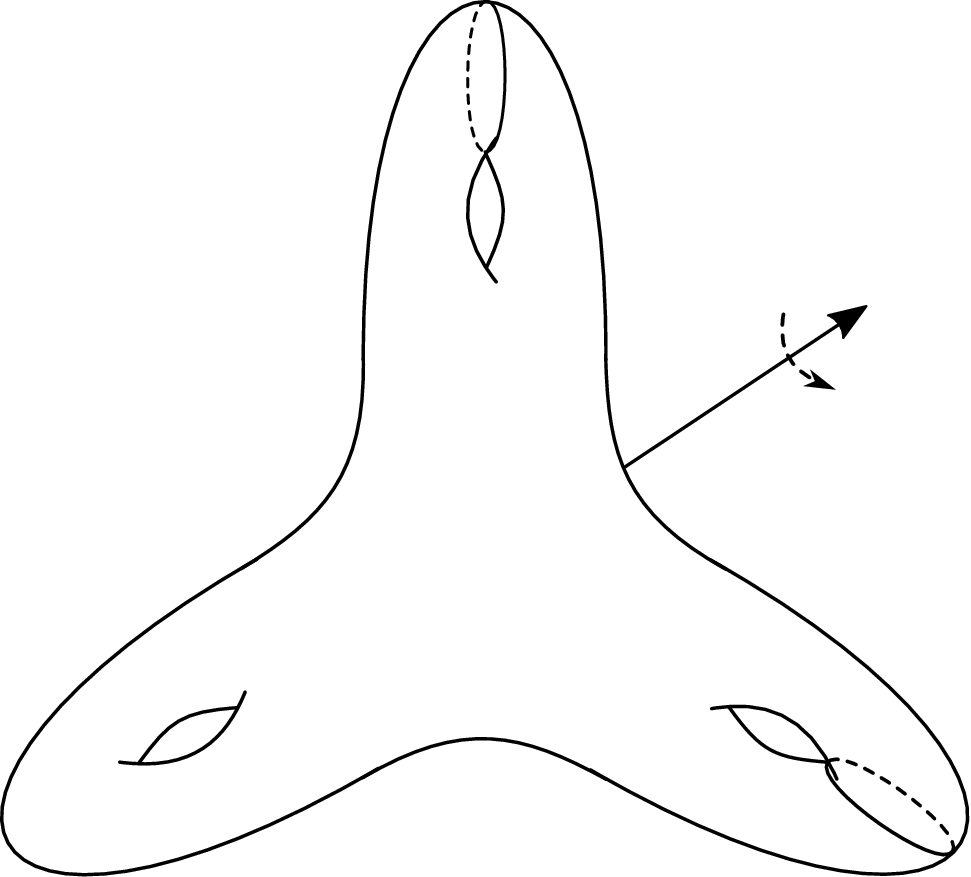}
\end{subfigure}
\caption{Some rotational cyclic actions on surfaces.}
\label{fig:rotational}
\end{figure}

Now, suppose that $F'$ is a non-rotational cyclic action. First, we consider the case when $D_{F'}$ is an irreducible Type 1 action. In~\cite[Theorem 2.7]{rajeevsarathy1}, it was shown that an irreducible Type 1 cyclic action $D$ can be realized as the $2\pi k/n(D)$-rotation (for some $k$ with $\gcd(k,n(D))=1$) of a canonical semi-regular hyperbolic polygon $2n(D)$-gon (or an $n(D)$-gon) with a suitable side-pairing. Thus, $D_{F'}$ can be realized as a multiple of $2\pi/n(D_{F'})$-rotation of a canonical hyperbolic polygon $P$. Let $c$ be the line segment joining the mid-points of two sides which are identified. Then $c$ is an essential simple closed curve in $S$. (We note that such a $c$ exists as $g(D_{F'})\neq 0$.) 

Since $D_{F'}$ is an irreducible Type 1 action, we have $n(D_{F'})>2$. Therefore, $c$ and $F'(c)$ are two distinct essential simple closed curves. We observe that $i(c,F'(c))=i(c,F(c))\leq 1$. If $c$ is separating, or $i(c,F'(c))=1$, or $i(c,F'^2(c))=1$, then the result follows from Lemma \ref{lem:wellsuited_curve} (we note that the normal closure of $F^2$ is contained in the normal closure of $F$). Assume that $c$ is nonseparating and $i(c,F'(c))=i(c,F'^2(c))=0$. We claim that $[c]$, $[F'(c)]$, and $[F'^2(c)]$ are not all equal. Since $i(c,F'(c))=i(c,F'^2(c))=0$, there is region of $P$ bounded by the three curve $c$, $F'(c)$, and $F'^2(c)$. But this is not possible as shown in Figure \ref{fig:rotation}. The result now follows from Lemma \ref{lem:wellsuited_curve}.
\begin{figure}[H]
\centering
\includegraphics[scale=.35]{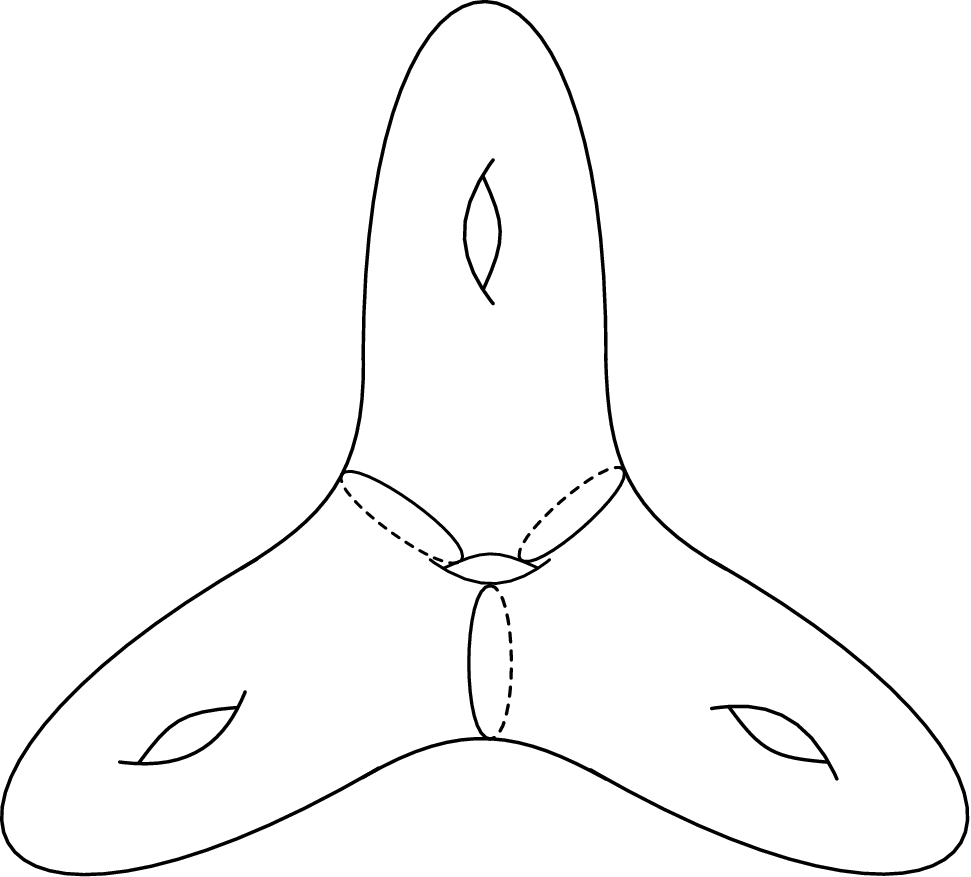}
\caption{Three disjoint homologous pairwise non-isotopic nonseparating essential simple closed curves in $S_4$.}
\label{fig:rotation}
\end{figure}

Now, we assume that $D_{F'}$ is of Type 2. It was proven in~\cite[Theorem 2.24]{rajeevsarathy1} that a Type 2 cyclic action $D$ can be constructed inductively from irreducible Type 1 actions by performing certain constructions that we termed as ``compatibilities". The key idea behind these compatibilities involved the deletion of disjoint (cyclically permuted) disks around a pair of orbits induced by an action (or a pair of actions) where the locally induced rotation angles add up to $0$ modulo $2\pi$. Consequently, $D_{F'}$ can be constructed from these compatibilities involving irreducible Type 1 periodics. Now, the result follows by considering any irreducible Type 1 component of $D_{F'}$.
\end{proof}
\end{prop}

It is not hard to verify that some well-studied classes of multitwists do not normally generate $\map(S_g)$. For example, given a simple closed curve $c$ and a bounding pair $\{a,b\}$ in $S_g$, the mapping classes $T_c^k$ when $|k|>1$ and $T_aT_b^{-1}$ do not generate $\map(S_g)$ normally. However, we will show that these mapping classes always have roots that normally generate $\map(S_g)$.

\begin{prop}
For $g\geq 3$, let $c$ be a simple closed curve, and let $\{a,b\}$ be a bounding pair in $S_g$. There are roots of $T_c^k$ for $k>1$ and $T_aT_b^{-1}$ that normally generate $\map(S_g)$.
\begin{proof}
In the case of a nonseparating curve $c$, the pseudo-periodic data set
\begin{center}
$\mathcal{D}_1=\llbracket 4(g-1);(D_1,1);-;(1,2(g-1),2(g-1);1,1)\rrbracket$
\end{center}
corresponds to a conjugacy class of a root of $T_c^{2(g-1)}$ of degree $4(g-1)$, where
\begin{center}
$D_1=(4(g-1),0;(1,4(g-1)),(2g-3,4(g-1)),(1,2)).$
\end{center}
Now we consider the case when $c$ is a separating curve in $S_g$ which separates the surface into two subsurfaces of positive genus $g_1$ and $g_2$. Consider the cyclic data sets
\begin{align*}
&D_1=(4g_1,0;(a_1,2),(a_1(2g_1-1),4g_1),(a_1,4g_1)_1) \text{ and}\\
&D_2=(4g_2+2,0;(a_2,2),(a_2g_2,2g_2+1),(a_2,4g_2+2)_1).
\end{align*}
Taking $n_1=4g_1$, $n_2=4g_2+2$, and $n=\lcm(n_1,n_2)$, since both $n_i$'s are not odd and $\gcd(n/n_1,n/n_2)=1$, it follows from Lemma \ref{lem:gcd_lemma} that there are integers $p_1$ and $p_2$ such that $\gcd(p_i,n_i)=1$ and
\begin{center}
$\frac{n}{n_1}p_1+\frac{n}{n_2}p_2\equiv k\pmod n$, where $\gcd(k,n)=1$.
\end{center}
For $a_i=p_i^{-1}\pmod {n_i}$, the fixed points of $D_i$ with the same subscript are compatible with twist factor $k$. Thus, the pseudo-periodic data set
\begin{center}
$\mathcal{D}_2=\llbracket n;(D_1,1),(D_2,1);\llparenthesis 1,k,k;1,2\rrparenthesis;-\rrbracket$
\end{center}
corresponds to a conjugacy class of a root of $T_c^k$ of degree $n$.

Now let $C=\{a,b\}$ be a bounding pair in $S_g$. In the proof of Proposition \ref{prop:torelli_levelm}, a root of bounding pair map $T_aT_b^{-1}$ of degree $g(g-2)$ has been constructed on surfaces of even genus $g$. Here we construct a root of a bounding pair map on the surfaces of odd genus. Let $s$ be a positive integer. When $g=3$, we consider the following data sets
\begin{align*}
&D_1=(4,1;(1,4)_1,(3,4)_2) \text{ and}\\
&D_2=(2,0;(1,2)_1,(1,2)_2,(1,2),(1,2)).
\end{align*}
When $g=4s+3$, for $g_1=2s+2$ and $g_2=2s$, we consider the following data sets
\begin{align*}
&D_1=(4g_1,0;(1,2),(1,4g_1)_1,(2g_1-1,4g_1)_2) \text{ and}\\
&D_2=(4g_2,0;(1,2),(4g_2-1,4g_2)_1,(2g_2+1,4g_2)_2).
\end{align*}
When $g=4s+1$ and $s$ is even, for $g_1=2s$ and $g_2=2s$, we consider the following data sets
\begin{align*}
&D_1=(g_1,1;(g_1-s+1,g_1)_1,(s-1,g_1)_2) \text{ and}\\
&D_2=(2,0;(1,2)_1,(1,2)_2,\underbrace{(1,2),(1,2),\dots,(1,2)}_{(2g_2)-\text{times}}).
\end{align*}
When $g=4s+1$ and $s$ is odd, for $g_1=2s+2$ and $g_2=2s-2$, we consider the following data sets
\begin{align*}
&D_1=(g_1,1;(1,g_1)_1,(g_1-1,g_1)_2) \text{ and}\\
&D_2=(g_2,1;(g_2-1,g_2)_1,(1,g_2)_2).
\end{align*}
In the above mentioned cases, the fixed points with the same subscripts are compatible with twist factor $\pm 1$, and hence, the pseudo-periodic data set
\begin{center}
$\mathcal{D}_3=\llbracket n;(D_1,1),(D_2,1);\llparenthesis 1,1,1;1,2\rrparenthesis,\llparenthesis 1,-1,-1;1,2\rrparenthesis;-\rrbracket$
\end{center}
corresponds to a root of $T_aT_b^{-1}$ of degree $n$, where $n=\lcm(n(D_1),n(D_2))$. From Proposition \ref{prop:pp_normal}, it follows that any pseudo-periodic mapping class conjugate to $\mathcal{D}_1$, $\mathcal{D}_2$, and $\mathcal{D}_3$ normally generate $\map(S_g)$. 
\end{proof}
\end{prop}

\section*{Acknowledgment} The authors would like to thanks Prof. Dan Margalit for inspiring the pursuit of this problem and also for sharing some helpful comments on the manuscript. The authors are also grateful to the referees for helpful suggestions that have significantly improved the exposition in this manuscript. The first author was supported by the Prime Minister Research Fellowship (PMRF) scheme instituted by the Ministry of Education, India.
\bibliographystyle{abbrv}
\bibliography{pseudo_periodic}
\end{document}